\documentclass[journal]{IAENGtran}
\ifCLASSINFOpdf
   \usepackage[pdftex]{graphicx}
   \DeclareGraphicsExtensions{.pdf,.jpeg,.png}
\else
   \usepackage[dvips]{graphicx}
   \DeclareGraphicsExtensions{.eps}
\fi

\usepackage{amssymb}
\usepackage{subfigure, graphicx}
\usepackage{latexsym}
\usepackage{graphicx}
\usepackage{amssymb, amsmath}
\usepackage{bm}

\newtheorem{lem}{Lemma}[section]
\newtheorem{vet}{Theorem}[section]

\newtheorem{tvr}{Proposition}[section]

\newcommand{\R}{\mathbb{R}}
\newcommand{\CC}{\mathbb{C}}
\newcommand{\Z}{\mathbb{Z}}
\newcommand{\N}{\mathbb{N}}

\newcommand{\K}{\mathcal{K}}

\def\vecx{\bm{x}}

\def\vecn{\bm{n}}
\def\vect{\bm{t}}

\newcommand{\mb}[1]{\mathbb{#1}}

\def\prf{\smallskip \noindent{\it Proof}.\ } 
\def\qed{{\hfill\QED \medskip}}
\renewcommand{\QED}{\hbox{\rule[0pt]{3pt}{6pt}}}

\def\ds{\hbox{d}s}
\def\du{\hbox{d}u}
\def\dnu{\hbox{d}\nu}

\begin{document}
%
\title{Application of the Enhanced Semidefinite Relaxation Method to Construction of the Optimal Anisotropy Function}

\author{Daniel \v{S}ev\v{c}ovi\v{c} and M\'aria Trnovsk\'a
\thanks{Manuscript received December 31, 2013; revised June 20, 2014. 
This work was supported in part by the APVV grant  SK-PT-0009-12 and 7FP EU STRIKE project No. 304617.}
\thanks{D. \v{S}ev\v{c}ovi\v{c} and M. Trnovsk\'a are with the 
Department of Applied Mathematics and Statistics, Comenius University, 842 48 Bratislava, Slovak Republic, e-mail: \{sevcovic,trnovska\}@fmph.uniba.sk}
}

\maketitle

\pagestyle{empty}
\thispagestyle{empty}

\begin{abstract}
In this paper we propose and apply the enhanced semidefinite relaxation technique for solving a class of non-convex quadratic optimization problems. The approach is based on enhancing the semidefinite relaxation methodology by complementing linear equality constraints by quadratic-linear constrains. We give sufficient conditions guaranteeing that the optimal values of the primal and enhanced semidefinite relaxed problems coincide. We apply this approach to the problem of resolving the optimal anisotropy function. The idea is to construct an  optimal anisotropy function as a minimizer for the anisotropic interface energy functional for a given Jordan curve in the plane. We present computational examples of resolving the optimal anisotropy function. The examples include boundaries of real snowflakes.
\end{abstract}

\begin{IAENGkeywords}
Enhanced semidefinite relaxation method, semidefinite programming, anisotropy function, Wulff shape
\end{IAENGkeywords}

\section{ Introduction}
\label{sec-1}

\IAENGPARstart{I}{n}
this paper we propose and apply the enhanced semidefinite relaxation technique  for solving the following nonlinear optimization problem:
\begin{equation}\label{P}
\begin{array}{rl}
\displaystyle\min & x^TP_0x+2q_0^Tx+r_0 \\
{\rm s.t.} & x^T P_l x+2q_l^Tx+r_l \leq 0, \quad l=1,\cdots, d, \\  
 & Ax=b, \\
 & H_0+\sum_{j=1}^n x_j H_j \succeq 0,
\end{array} 
\end{equation}
where $x\in\mb{R}^n$ is the variable and the data: $P_0, P_l$ are $n\times n$ real symmetric matrices, $q_0, q_l\in\mb{R}^n$, $r_0, r_l\in\mb{R}$, $A$ is an $m\times n$ real matrix of the full rank, $b\in\mb{R}^m$ and $H_0, H_1, \cdots, H_n$ are $k\times k$ complex Hermitian matrices. The last constraint in \eqref{P} is referred to as the linear matrix inequality (LMI). The relation $H\succeq 0$ means that the matrix $H$ is positive semidefinite. The optimal value of problem \eqref{P} will be denoted by $\hat p_1$. Optimization problems of the form \eqref{P} arise from various applications of combinatorial optimization, engineering, physics and other fields of applied research. In these problems the objective function $x\mapsto  x^TP_0x+2q_0^Tx+r_0$ need not be necessarily convex, in general. 

Our aim is to propose and then apply a novel method for solving non-convex optimization problems of the form \eqref{P}. The method is based on enhancing the usual semidefinite relaxation methodology by complementing linear equality constraints by quadratic-linear constraints. In the usual semidefinite relaxation procedure the quadratic term $X:=x x^T$ is relaxed by the linear matrix inequality $X\succeq x x^T$.  Notice that for the case there are neither linear constraints ($A=0$) nor LMI constraints ($H_j=0, j=0,1,\cdots, n$) and $d=1$ the method of semidefinite relaxation for \eqref{P} was analyzed in \cite[Appendix C.3]{bova}. We also refer the reader to papers by Boyd and Vanderberghe \cite{bovasdp,sdprelax}, Bao et al. \cite{review}, Nowak \cite{nowak} and Shor \cite{shor} for an overview of semidefinite relaxation techniques for solving various classes of non-convex quadratic optimization problems. Optimal solutions of a second order cone programming problem with box and linear constraints of have been studied by Hasuike in \cite{H2011}. 

Our idea of enhancing such a semidefinite relaxation technique consists in adding the linear constraint $AX= b x^T$ between the unknown vector $x$ and the semidefinite relaxation $X\succeq x x^T$.   Let us emphasize that usual semidefinite relaxation techniques just replace the quadratic terms forming the matrix $x x^T$ by an $n\times n$ matrix $X$ such that $X\succeq x x^T$. In our contribution, we propose to enhance such a relaxation by adding a new constraint $A X = b x^T$ which can be deduced from $Ax=b$ in the case $X=x x^T$. With regard to Proposition~\ref{tight} this constraint binds matrices $X$ and $x x^T$ to be close to each other in the sense of the rank-defect of their difference.

We apply this method to the problem of construction the optimal anisotropy function. The anisotropy function $\sigma$ describing the so-called Finsler metric in the plane occurs in  various models from mathematical physics. It particular, it enters the anisotropic Ginzburg-Landau free energy and the nonlinear parabolic Allen-Cahn equation with a diffusion coefficient depending on $\sigma$ (cf. Belletini and Paolini \cite{belletini}, Bene\v{s} et al. \cite{benes2003,hong2013}). It is also important in the field of differential geometry and its applications to anisotropic  motion of planar interfaces (see e.g. Perona and Malik \cite{perona}, Weickert \cite{weickert}, Mikula and the first author \cite{MS2004}). In all aforementioned models the anisotropy function $\sigma$ enters the model as an external given data. On the other hand, considerably less attention is put on understanding and construction the anisotropy function $\sigma$ itself. In the application part of this paper we present a novel idea how to construct an optimal anisotropy function by means of minimizing the total anisotropic interface energy $L_\sigma(\Gamma)$ for a given Jordan curve $\Gamma$ in the plane. It leads to a solution to the optimization problem: 
$\inf_\sigma L_\sigma(\Gamma)$ where $L_\sigma(\Gamma) = \int_\Gamma \sigma(\nu)\ds$. 
Here the unknown anisotropy function $\sigma$ is a nonnegative function of the tangent angle $\nu$ of the curve $\Gamma$. The tangent angle $\nu$ is defined by the relation: $\vecn=(-\sin\nu, \cos\nu)^T$, where $\vecn$ is the unit inward normal vector to $\Gamma$. In this paper we show how such an optimization problem can be reformulated as a non-convex quadratic programming problem  with linear matrix inequalities of the form \eqref{P}. The method of the enhanced semidefinite relaxation of \eqref{P} can be also used in other applications leading to non-convex constrained problems.

The paper is organized as follows. In Section II we present the enhanced semidefinite relaxation method for solving the optimization problem \eqref{P}. The method is based on enhancing the classical semidefinite relaxation methodology by means of complementation of linear equality constraints by quadratic-linear constrains. We give sufficient conditions guaranteeing that the optimal values of primal and enhanced semidefinite relaxed problems coincide. In Section III we investigate the problem of construction of the optimal anisotropy function minimizing the total anisotropic interface energy of a given Jordan curve in the plane. We propose two different criteria for the anisotropy function based on the linear and second order quadratic type of constraints. Section IV is devoted to representation of the optimal anisotropy problem by means of the Fourier series expansion of the anisotropy function. In Section V we show that the optimization problem is semidefinite representable and it fits into the general framework of the class of non-convex optimization problems having the form \eqref{P}. In Section VI we present several numerical experiments for construction of the optimal anisotropy function for various Jordan curves in the plane including, in particular, boundaries of real snowflakes.

\section{Enhanced semidefinite relaxation method}

Our aim is to investigate problem \eqref{P} by means of methods and techniques of non-convex optimization. Namely, we will apply a method of the semidefinite relaxation of \eqref{P} in combination with complementation of \eqref{P} by quadratic-linear constraints. The method will be referred to as the {\it enhanced semidefinite relaxation method} for solving \eqref{P}. In the recent paper \cite{ST} the theoretical and numerical aspects of the method have been investigated.

For a transpose of the matrix $X$ we will henceforth write $X^T$. For complex conjugate of a complex matrix $H$ we will write $H^*$, i.~e. $H^*=\bar{H}^T$. The sets of real $n\times n$ symmetric ($A=A^T$) and complex Hermitian matrices ($H=H^*$) are denoted by $\mathcal{S}^n$ and $\mathcal{H}^n$, respectively. We will write $A\succeq 0$ ($A\succ 0$), if a real symmetric matrix or a complex  Hermitian matrix $A$ is positive semidefinite (positive definite).

First, it should be obvious that  problem \eqref{P} is equivalent to the following augmented problem:
\begin{equation}\label{P1}
\begin{array}{rl}
\displaystyle\min & x^TP_0x+2q_0^Tx+r_0 \\
\hbox{s. t.} & x^T P_l x+2q_l^Tx+r_l \leq 0, \quad l=1,\cdots, d, \\  
 & Ax=b, \ \  Ax x^T=bx^T,\\
 & H_0+\sum_{j=1}^n x_j H_j \succeq 0. 
\end{array} 
\end{equation}
It contains the additional quadratic-linear constraint $Ax x^T=bx^T$ which clearly follows from the equality constraint $Ax=b$.  The optimal value of \eqref{P1} is equal to the value $\hat p_1$. The method of semidefinite relaxation is based on the idea that all the terms in \eqref{P1} of the form $X=xx^T$ are relaxed by a convex constraint $X\succeq x x^T$, i.~e. the matrix  $X-x x^T$ is positive semidefinite. Enhancement of \eqref{P} means that the equality constraint $Ax=b$ will be complemented by the quadratic-linear constraint $Ax x^T = b x^T$. This is a dependent constraint. On the other hand, in combination with the semidefinite relaxation of $X=x x^T$, the original quadratic-linear constraint $Ax x^T=b x^T$ will be transformed to the set of linear equations $AX = b x^T$ between the unknown vector $x$ and the $n\times n$ semidefinite relaxation matrix $X$ such that $X\succeq xx^T$. 

Since $x^T P_l x = \hbox{tr}(x^T P_l x) =  \hbox{tr}(P_l x x^T)$ construction of a  semidefinite relaxation of \eqref{P1} is rather simple and it consists in relaxing the equality $X=x x^T$ by the semidefinite inequality $X\succeq x x^T$. The enhanced semidefinite relaxation of \eqref{P} now reads as follows:
\begin{equation}\label{D3}
\begin{array}{rl}
 \displaystyle\min & \hbox{tr}(P_0 X)+2q_0^Tx+r_0\\
\hbox{s. t. } & \hbox{tr}(P_l X)+2q_l^Tx+r_l\leq 0, \ \ l=1,\cdots,d,\\
    & Ax=b, \ AX=bx^T, \  X\succeq x x^T,\\
& {H}_0+\sum_{j=1}^n x_j {H}_j\succeq 0.
\end{array}
\end{equation}
Notice that, using the property of the Schur complement the inequality $X\succeq x x^T$ can be rewritten as the linear matrix inequality, i.~e.
\[
X\succeq x x^T  \Longleftrightarrow  
\left(\begin{matrix}
 X   & x \\ 
 x^T & 1 
\end{matrix}\right) \succeq 0
\]
(cf. Zhang \cite{Z}). The optimal value of \eqref{D3} will be denoted by $\hat p_2$. Since the resulting semidefinite relaxed problem \eqref{D3} is a convex optimization problem it can be efficiently solved by using available solvers for nonlinear programming problems over symmetric cones, e.~g. SeDuMi or SDPT3 solvers \cite{sturm}. 

\subsection{Equivalence of problems \eqref{P} and \eqref{D3}}

In this section we provide a sufficient condition guaranteeing that the primal problem \eqref{P1} and its enhanced semidefinite relaxation \eqref{D3} yield the same optimal values $\hat p_1$ and $\hat p_2$, respectively. First we compare optimal values of \eqref{P} and  \eqref{D3}. Next, under additional assumptions, we show equivalence of optimal values and optimal solutions to \eqref{P} and \eqref{D3}.

\begin{vet}\label{tvr1}
Suppose that problem \eqref{P} is feasible. Then the enhanced semidefinite relaxation problem \eqref{D3} is also feasible. For the optimal values $\hat p_1$ of \eqref{P} and $\hat p_2$ of \eqref{D3} we have $\hat p_1\geq \hat p_2$.
\end{vet}

\prf
Let $x\in\R^n$ be a feasible solution to \eqref{P}. Clearly, $x$ is feasible to the augmented problem \eqref{P1} as well. Set  $X=x x^T$. As $\hbox{tr}(P_l X) = \hbox{tr}(P_l x x^T) = \hbox{tr}(x^T P_l x) = x^T P_l x$ for $l=0,1, \dots, d$, we have that the pair $(x,X)$ is feasible to \eqref{D3}. Finally, 
$x^T P_0 x + 2q_0^T x + r_0 = \hbox{tr}(P_0 X) + 2q_0^T x + r_0 \ge \hat p_2$ because $(x,X)$ is feasible to \eqref{D3}. Hence $\hat p_1\geq \hat p_2\ge -\infty$. \qed

\begin{vet}\label{prop-ekviv}
Assume $P_l\succeq 0, l=1,\dots,d$, are positive semidefinite matrices. Suppose that $(\hat x, \hat X)$ is an optimal solution to \eqref{D3} satisfying the inequality 
\begin{equation}
\hbox{tr}(P_0 \hat X) \ge \hat x^T P_0 \hat x. 
\label{gap}
\end{equation}
Then $\hat x$ is an optimal solution to \eqref{P} and $\hat p_1 = \hat p_2$. Moreover, we have  $\hbox{tr}(P_0 \hat X) = \hat x^T P_0 \hat x$.
\end{vet}

\prf
It follows from basic properties of positive semidefinite matrices that 
\begin{equation}
\hbox{tr}(P M) \ge 0 \quad\hbox{for any} \ P\succeq 0, M\succeq 0, \ P,M\in {\mathcal S}^{n}.
\label{ineq}
\end{equation}
Now, if $(\hat x, \hat X)$ is an optimal solution to \eqref{D3} then $M= \hat X - \hat x \hat x^T \succeq 0$. Then  $\hbox{tr}(P_l (\hat X- \hat x \hat x^T) ) \ge 0$ and so $\hat x^T P_l \hat x \le  \hbox{tr}(P_l \hat X)$ for $l=1,\dots,d$. Hence $\hat x$ is a feasible solution to \eqref{P}. Taking into account inequality \eqref{gap} and Theorem~\ref{tvr1} we conclude $\hat p_2 = \hbox{tr}(P_0 \hat X) + 2q_0^T \hat x + r_0  \ge \hat x^T P_0 \hat x + 2q_0^T \hat x + r_0 \ge \hat p_1 \ge \hat p_2$. Therefore $\hat p_1 = \hat p_2$ and $\hbox{tr}(P_0 \hat X) = \hat x^T P_0 \hat x$, as claimed. 
\qed

In the next proposition we give a sufficient condition guaranteeing inequality \eqref{ineq}. It is closely related to the Finsler characterization of positive semidefinitness of a matrix $P_0$ over the null subspace $\{ x\in\R^n | Ax=0\}$ (cf. \cite{finsler}).

\begin{tvr}\label{Finsler}
The inequality $\hbox{tr}(P_0  X) \ge  x^T P_0  x$ is satisfied by any $(x,X)$ feasible to \eqref{D3} provided that there exists $\varrho\in \R$ such that 
$P_0 + \varrho\, A^T A \succeq 0$.
\end{tvr}

\prf
Suppose that $P_0 + \varrho\, A^T A \succeq 0$ for some $\varrho$. By \eqref{ineq} we have 
$0 \le \hbox{tr}( (P_0 +\varrho A^T A) (X-x x^T) ) 
= \hbox{tr}(P_0 X - P_0 x x^T)  + \varrho \, \hbox{tr}( A^T [ AX - Ax x^T]) 
=  \hbox{tr}(P_0 X)  - x^T P_0 x$,
because $AX = b x^T = Ax x^T$ for any $(x,X)$ feasible to \eqref{D3}.
\qed

Finally, we show that any feasible solution $(x,X)$ to the enhanced semidefinite relaxation problem \eqref{D3} is tight in the sense that the gap matrix $X-x x^T$ is a positive semidefinite matrix of the rank at most of $n-m$.

\begin{tvr}\label{tight}
Suppose that an $m\times n$ real matrix $A$ has the full rank $m$. Then $rank(X-x x^T) \le n-m$ for any feasible solution $(x,X)$ to \eqref{D3}. In particular, $X\succeq x x^T$ but $X\not\succ x x^T$. 
\end{tvr}

\prf Let $Y:= X -x x^T$. Then $Y\succeq 0$ and 
$A Y = AX - Ax x^T = (b-Ax) x^T =0$. Since $AY=0$ the range $S(Y)$ of the matrix $Y$ is a subspace of the null space $N(A)$ of the matrix $A$. Thus $\hbox{rank\,} (Y) =\hbox{dim\,} S(Y) \le \hbox{dim\,} N(A) = n-m$,
as claimed. \qed

\section{Application of the method for construction of optimal anisotropy function}

In many applications arising from material science, differential geometry, image processing knowledge of the so-called anisotropy function $\sigma$ plays an essential role. In the case of the Finsler geometry of the plane the anisotropy function $\sigma=\sigma(\nu)$ depends on the tangent angle $\nu$ of a curvilinear boundary $\Gamma$ enclosing a two dimensional connected area. The total anisotropic interface energy $L_\sigma(\Gamma)$ of a closed curve $\Gamma \in \R^2$ can be defined as follows
$L_\sigma(\Gamma) = \int_\Gamma \sigma(\nu) \ds$.
The anisotropy function $\sigma$ is closely related to the fundamental notion describing the generalized (Finsler) geometry in the plane. Such a geometry can be  characterized by the so-called Wulff shape $W_\sigma$. Given a $2\pi$-periodic nonnegative anisotropy function $\sigma=\sigma(\nu)$ the Wulff shape in the plane is defined as $W_{\sigma}=\bigcap_{\nu\in[0,2\pi]}\left\{\vecx \ |  \ -\vecx^T \vecn \leq \sigma(\nu)\right\}$, 
where $\vecn=(-\sin\nu, \cos\nu)^T$ is the unit inward vector. It is well known that the boundary $\partial W_\sigma$ can be parameterized as follows: $\partial W_{\sigma}=\left\{\vecx(\nu)\ |\ \vecx(\nu)=-\sigma(\nu)\vecn+\sigma'(\nu)\vect, \ \nu\in  [0,2\pi]\right\}$ 
where $\vect\equiv(t_1,t_2)^T=(\cos\nu, \sin\nu)^T$ is the unit tangent vector to the boundary $\partial W_\sigma$ of the Wulff shape. 
Its curvature $\kappa$ is given by $\kappa=[\sigma(\nu)+\sigma''(\nu)]^{-1}$ (see \cite{SY} for details). Hence the Wulff shape $W_{\sigma}$  is a convex set if and only if $\sigma(\nu)+\sigma''(\nu) \ge 0$ for all $\nu\in \R$. Henceforth,  we will assume the anisotropy function $\sigma\in\K$ belongs to the cone of $2\pi$-periodic functions
\begin{equation}
\K=\{\sigma\in W^{2,2}_{per}(0,2\pi)\ |\ \sigma\ge 0, \ \sigma + \sigma'' \ge 0\},
\label{cone}
\end{equation} 
where  $W^{2,2}_{per}(0,2\pi)$ denotes the Sobolev space of all real valued $2\pi$-periodic functions having their distributional derivatives square integrable up to the second order.

For a boundary $\partial W_\sigma$ of the convex set $W_\sigma$ the tangent angle $\nu$ can be used as a parameterization of $\partial W_\sigma$. Moreover, as $\kappa=\partial_s \nu$ we can calculate the area $|W_\sigma|$ of the Wulff shape as follows:
\begin{eqnarray}
|W_\sigma| &=& -\frac{1}{2}\int_{\partial W_{\sigma}}\vecx^T \vecn\,\ds
=\frac{1}{2}\int_{\partial W_{\sigma}}\sigma(\nu)\,\ds
\nonumber \\
&= &\frac12 \int_0^{2\pi} \sigma(\nu) [\sigma(\nu) + \sigma''(\nu)] \hbox{d}\nu 
\label{areaW}\\
&=& \frac12 \int_0^{2\pi} |\sigma(\nu)|^2 - |\sigma'(\nu)|^2 \hbox{d}\nu,
\nonumber
\end{eqnarray}
because $\dnu = \kappa \ds = [\sigma + \sigma'']^{-1} \ds$. If $\sigma\equiv 1$ then the boundary $\partial W_{1}$ of $W_1$ is a circle with the radius 1, and $|W_{1}|=\pi$.

The main contribution of this part of the paper is to propose a method how to construct the anisotropy function $\sigma$ with respect to minimization of the total interface energy $L_\sigma(\Gamma)$ provided that the Jordan curve (a closed $C^1$ smooth non-selfintersecting curve in the plane) is given. In practical applications, such a curve $\Gamma$ can represent a boundary of an important object like e.~g. a boundary of a snowflake for which anisotropic growth model we want to construct the underlying anisotropy function $\sigma$. Or, it may represent a boundary of a typical object in the image we want to segment by means of the anisotropic diffusion image segmentation model. 

In what follows, we will analyze two different approaches for construction of the optimal anisotropy function $\sigma$. We will show that imposing the first order constraint on the anisotropy function does not lead  to satisfactory results and the second order constraint should be taken into account when resolving the optimal anisotropy function $\sigma$.

First we notice the following homogeneity properties of the interface energy and the area of the Wulff shape hold true:
\begin{equation}
L_{t\sigma}(\Gamma) = t L_{\sigma}(\Gamma), \quad |W_{t\sigma}| = t^2 |W_\sigma|, \label{homog}
\end{equation}
for any $\sigma\in \K$ and all $t>0$. Moreover, $\inf_{\sigma\in\K} L_\sigma(\Gamma)=0$ for $\sigma\equiv 0\in \K$. In order to obtain a nontrivial anisotropy function minimizing the interface energy $L_\sigma(\Gamma)$ of a given curve $\Gamma$ we have to impose additional constraints on $\sigma$. We will distinguish two cases - linear and quadratic constraints on $\sigma$.  More precisely, given a Jordan curve $\Gamma$ in the plane we construct the optimal anisotropy function $\sigma$ as follows:

\begin{enumerate}
\item ({\it First order linear constraint imposed on $\sigma$}) 
\\ The anisotropy function $\sigma$ is a minimizer of 
\begin{equation}
\begin{array}{rl}
\displaystyle \inf_{\sigma\in \K} & L_\sigma(\Gamma) \\
{\rm s. t.} & \sigma_{avg} = 1,
\end{array} 
\label{first}
\end{equation}
where $\sigma_{avg}=\frac{1}{2\pi}\int_0^{2\pi}\sigma(\nu)\dnu$ is the average of $\sigma$.

\item ({\it Second order constraint imposed on $\sigma$}) 
\\ The anisotropy function $\sigma$ is a minimizer of 
\begin{equation}
\begin{array}{rl}
\displaystyle \inf_{\sigma\in\K} & L_\sigma(\Gamma) \\
{\rm s. t.} & |W_\sigma| = 1,
\end{array} 
\label{second}
\end{equation}
where $|W_\sigma|$ is the area of the Wulff shape.
\end{enumerate}

With regard to the homogeneity properties \eqref{homog}, the constrained problem \eqref{second}  can be also viewed as a solution to the inverse Wulff problem stated as follows: 
\[
\inf_{\sigma\in \K} \Pi_\sigma(\Gamma), \quad\hbox{where}\ \ \Pi_\sigma(\Gamma) = 
\frac{L_\sigma(\Gamma)^2}{4| W_\sigma| {\mathcal A}(\Gamma)}
\]
is the anisoperimetric ratio of a curve $\Gamma$ for the underlying anisotropy function $\sigma$. In \cite{SY} Yazaki and the author showed the following anisoperimetric inequality:
\begin{equation}
\frac{L_\sigma(\Gamma)^2}{4| W_\sigma| {\mathcal A}(\Gamma)} \ge 1,
\label{mixedanisoperim}
\end{equation}
where ${\mathcal A}(\Gamma)$ is the area enclosed by $\Gamma$. The equality is attained if and only if $\Gamma$ is homothetically similar to $\partial W_\sigma$. It is a generalization of the anisoperimetric inequality due to Wulff \cite{Wulff1901} (see also Dacorogna and Pfister \cite{dacorogna}) originally shown for $\pi$- periodic anisotropy function $\sigma$ only.

\section{Fourier series representation of the anisotropy function}

Since the anisotropy function $\sigma\in\K$ is a $2\pi$-periodic real function of a real variable $\nu\in\R$ it is useful to represent $\sigma$ by means of coefficients of its Fourier series expansion. A function $\sigma\in W^{2,2}_{per}(0,2\pi)$ can be represented by its complex Fourier series: 
\begin{equation}
\sigma(\nu) = \sum_{k=-\infty}^\infty \sigma_k e^{i k\nu}, \ \ 
\sigma_k=\frac{1}{2\pi}\int_0^{2\pi} e^{-i k\nu} \sigma(\nu)\hbox{d}\nu
\label{fourier}
\end{equation}
are complex Fourier coefficients. Since $\sigma(\nu)$ is assumed to be a real valued  function we have $\sigma_{-k} = \bar\sigma_k$ for any $k\in\Z$ and $\sigma_0\in\R$. 

In what follows, we will express the anisotropic interface energy $L_\sigma(\Gamma)$, the average value $\sigma_{avg}$ as well as the area $|W_\sigma|$ of the Wulff shape in terms of the Fourier coefficients $\sigma_k, k\in \Z$. Furthermore, we will provide a necessary and sufficient semidefinite representable condition for $\sigma$ to belong to the cone $\K$.

\subsection{Representation of the interface energy}

In terms of Fourier coefficients $\sigma_k, k\in \Z,$ the interface energy $L_\sigma(\Gamma)$ can be represented  as follows:
\begin{eqnarray}\label{interface-energy}
L_\sigma(\Gamma) &=& \int_\Gamma \sigma(\nu)\ds = \sum_{k=-\infty}^\infty \sigma_k \int_\Gamma e^{i k \nu}\ds \nonumber \\ 
&=& \sum_{k=-\infty}^\infty \bar c_k \sigma_k
= c_0 \sigma_0 + 2\Re  \sum_{k=1}^\infty \bar c_k \sigma_k, 
\end{eqnarray}
where the complex coefficients 
\begin{equation}
c_k:=\int_\Gamma e^{-i k \nu}\ds, \quad k\in\Z, 
\label{fouriercoeff}
\end{equation}
depend on the Jordan curve $\Gamma$ only. Using  the unit tangent vector $\vect=(t_1,t_2)^T=(\cos\nu, \sin\nu)^T$ to $\Gamma$ the coefficients $c_k, k\in \Z$, can be calculated as follows:
\[
c_k = \int_\Gamma e^{-ik\nu}\ds = \int_\Gamma (t_1 - i t_2)^k \ds\,.
\]
Notice that $c_0=\int_\Gamma \ds$ is the length $L(\Gamma)$ of the curve $\Gamma$. 

\begin{figure}
\begin{center}
\subfigure[]{\includegraphics[width=0.23\textwidth]{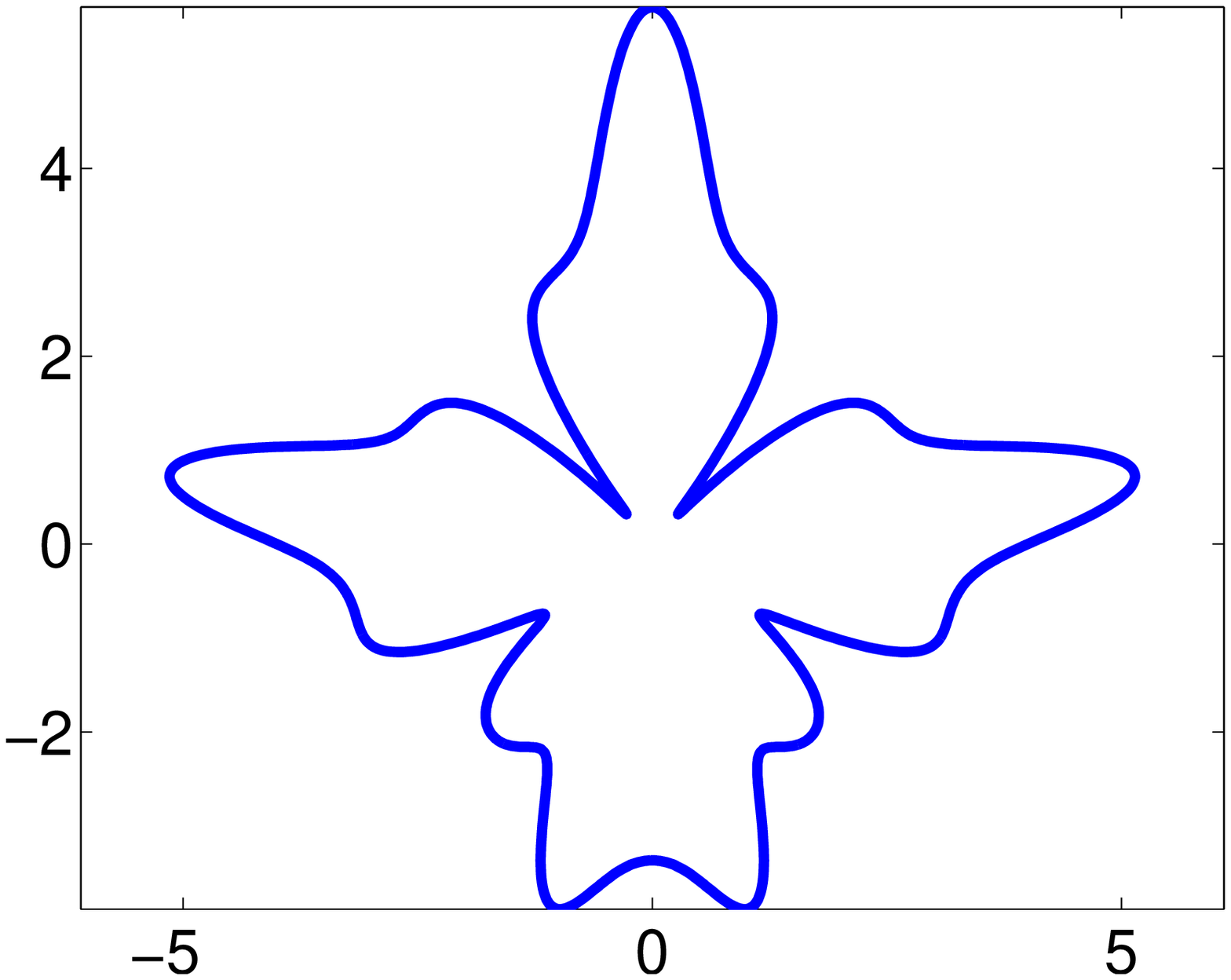}}
\subfigure[]{\includegraphics[width=0.25\textwidth]{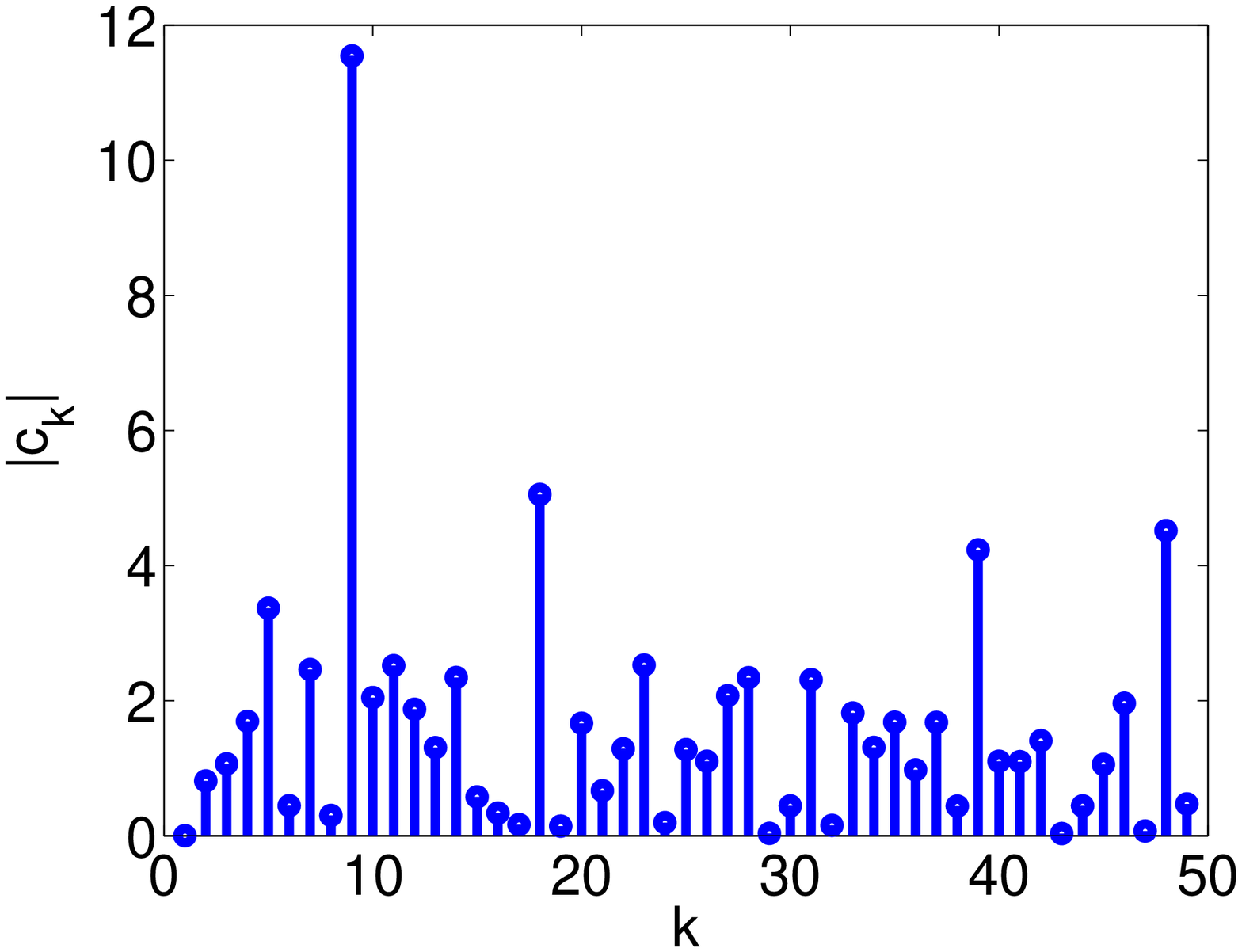}}
\\
\subfigure[]{\includegraphics[width=0.23\textwidth]{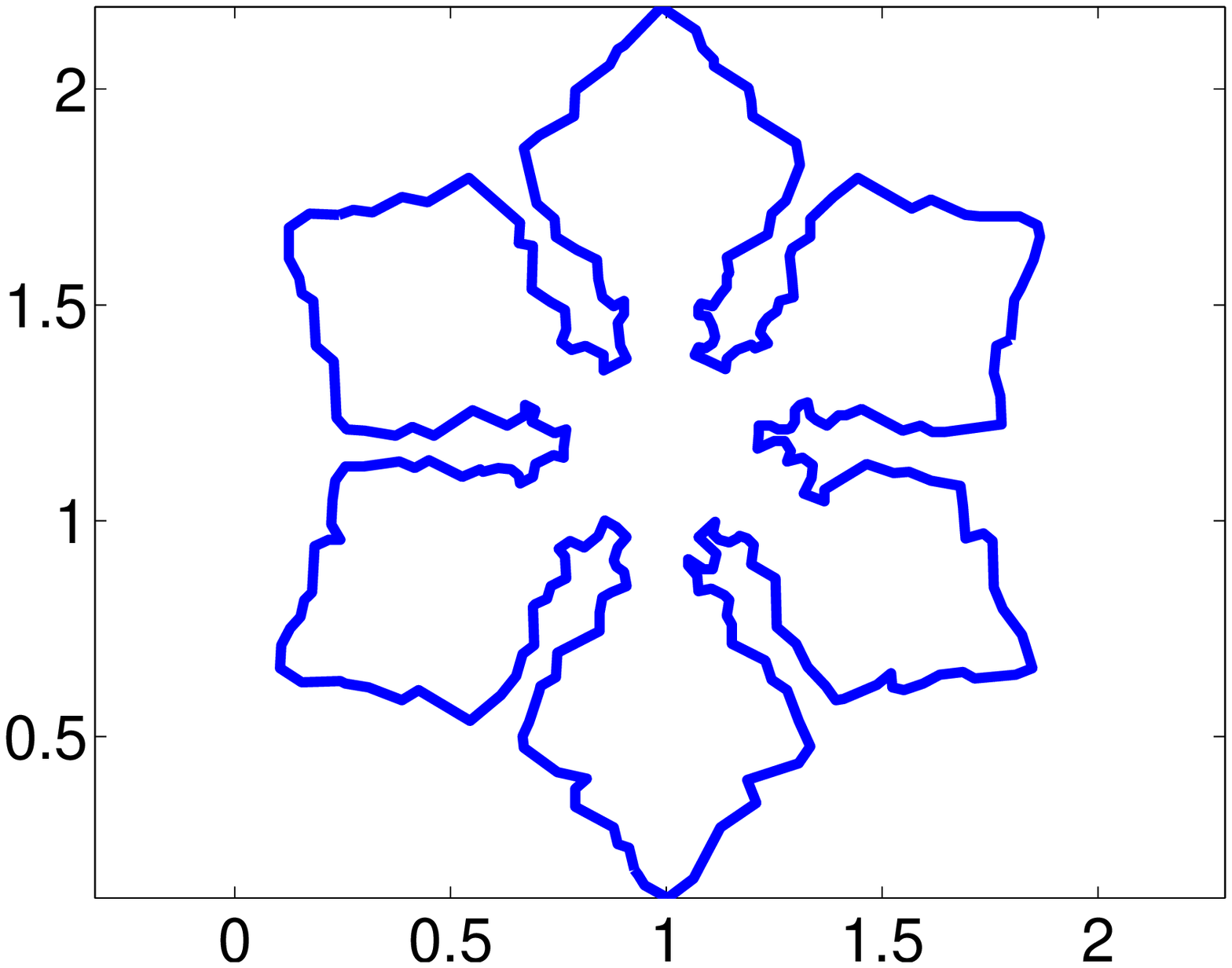}}
\subfigure[]{\includegraphics[width=0.25\textwidth]{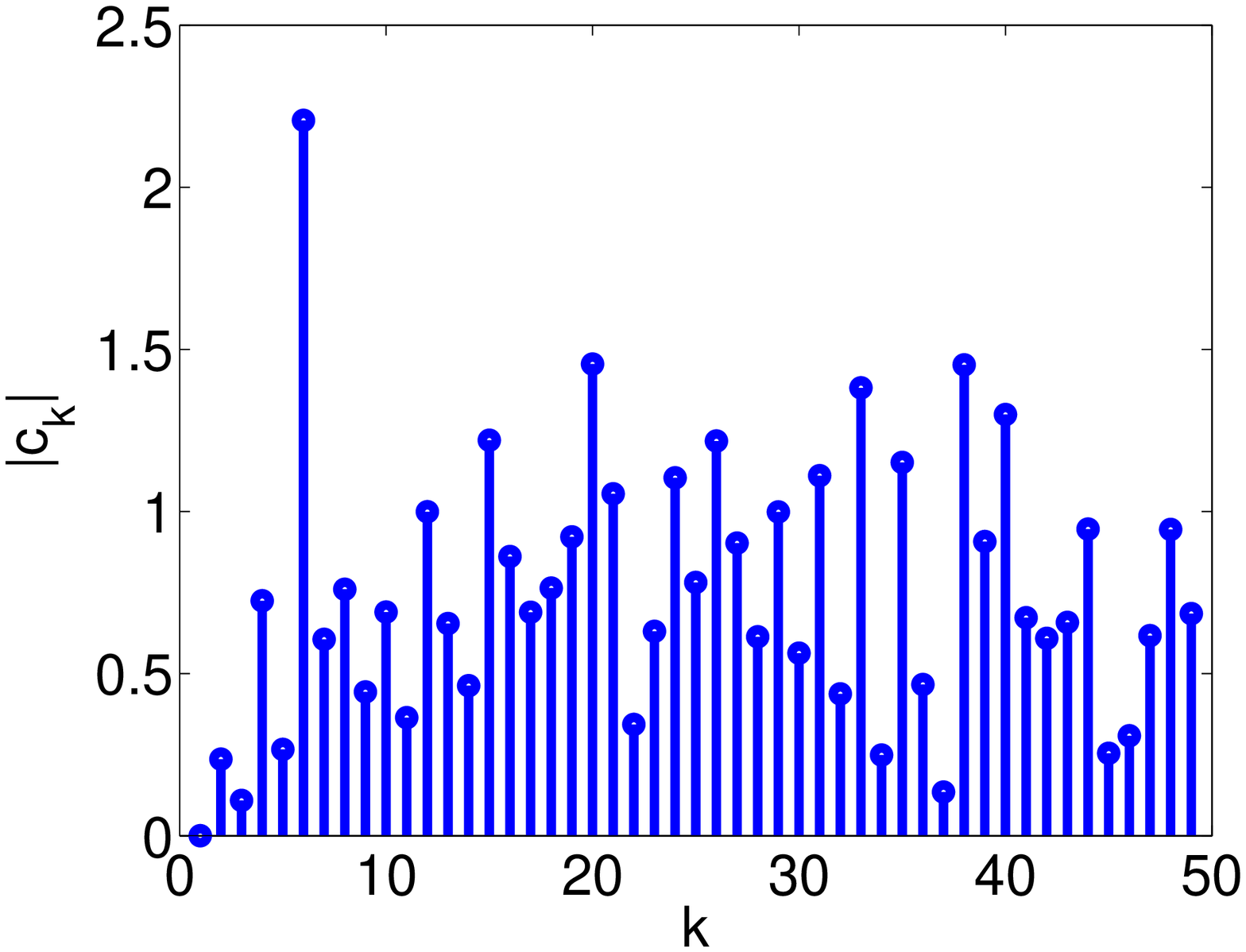}}
\end{center}
\caption{\small
Jordan curves (a,c) corresponding to the boundary of a dendrite and a real snowflake, respectively. The spectrum of moduli $|c_k|, k\ge 1,$ is shown in (c,d), respectively
}
\label{fig:fspec}
\end{figure}

In Fig.~\ref{fig:fspec} we plot moduli of $|c_k|, k\ge 1,$ of a dendrite type of a curve $\Gamma$ (top) and the boundary of a real snowflake (bottom). For the analytic description of the curve shown in Fig.~\ref{fig:fspec} (a), we refer to  Section VI.  In order to compute the above path integral, the curve $\Gamma$ was approximated by a polygonal curve $\hbox{poly}(\vecx^{(0)}, \vecx^{(1)}, \cdots \vecx^{(K)})$ with vertices $\vecx^{(0)}, \vecx^{(1)}, \cdots \vecx^{(K)}$. The unit tangent $\vect^{(j)}$ vector at $\vecx^{(j)}$ has been approximated by $\vect^{(j)} \equiv (\vecx^{(j+1)}-\vecx^{(j-1)})/\Vert \vecx^{(j+1)}-\vecx^{(j-1)}\Vert$. Here $\Vert \vecx\Vert$ is the Euclidean norm of a vector $\vecx$. Since $\ds =\Vert\partial_u \vecx\Vert \du \approx \frac12 \Vert \vecx^{(j+1)}-\vecx^{(j-1)}\Vert$ the coefficients $\{c_k, k\in \Z\}$ were approximated as follows:
\begin{eqnarray}\label{cp-discr}
c_k &=& \int_\Gamma (t_1 - i  t_2\bigl)^k\ds\nonumber  \\
&\approx& 
\frac12 \sum_{j=1}^{K-1} \bigr(t^{(j)}_1 - i  t^{(j)}_2\bigl)^k \Vert \vecx^{(j+1)}-\vecx^{(j-1)}\Vert.
\end{eqnarray}

\subsection{Representation of the average value of the anisotropy function}
The representation of the average value $\sigma_{avg}=\frac{1}{2\pi}\int_0^{2\pi}\sigma(\nu)\dnu$ is rather simple because 
\begin{equation}
\sigma_{avg}=\frac{1}{2\pi}\int_0^{2\pi}\sigma(\nu)\dnu 
=  \frac{1}{2\pi} \int_0^{2\pi} \sum_{k=-\infty}^\infty \sigma_k e^{i k\nu}= 
\sigma_0. 
\label{sigmaavg}
\end{equation}

\subsection{Representation of the Wulff shape area}

The area $|W_\sigma|$ of the Wulff shape can be easily expressed in terms of Fourier coefficients as follows:
\begin{eqnarray}
|W_\sigma| &=& \frac12 \int_0^{2\pi} |\sigma(\nu)|^2 - |\sigma'(\nu)|^2 \hbox{d}\nu \nonumber 
\\
&=& \frac12 \int_0^{2\pi} \sum_{k,m=-\infty}^\infty \bar\sigma_m \sigma_k (1- m k) e^{i(k-m)\nu} \dnu
\nonumber\\
&=& \pi \sum_{k=-\infty}^\infty (1-k^2) |\sigma_k|^2  \label{area-Wulff}
\\
&=& \pi\sigma_0^2 + 2\pi \sum_{k=1}^\infty (1-k^2) |\sigma_k|^2. \nonumber 
\end{eqnarray}

\subsection{Finite Fourier modes approximation}

In order to compute the optimal anisotropy function $\sigma$ we approximate $\sigma$ by its finite Fourier modes approximation up to the order $N$. To this end, we introduce the finite dimensional sub-cone $\K^N$ of $\K$ where
\begin{eqnarray}\label{coneN}
\K^N &=&\{ \sigma \in \K\ |\  \ \exists (\sigma_0, \sigma_1, \cdots,\sigma_{N-1})^T\in\CC^N,\nonumber \\
&& \sigma(\nu) = \sum_{k=-N+1}^{N-1} \sigma_k e^{i k \nu}\}.
\end{eqnarray}
Here $\sigma_{-k}=\bar\sigma_k$. For any $\sigma\in\K^N$ we have 
\begin{eqnarray*}
L_{\sigma}(\Gamma) &=& c_0 \sigma_0 + 2\Re  \sum_{k=1}^{N-1} \bar c_k \sigma_k,
\\
|W_{\sigma}| &=& \pi\sigma_0^2 + 2\pi \sum_{k=1}^{N-1} (1-k^2) |\sigma_k|^2.
\end{eqnarray*}

\subsection{Criteria for non-negativity of partial Fourier series}

Following the classical Riesz-Fejer factorization theorem (cf. \cite[pp. 117--118]{riesz}), in \cite{mclean} McLean and Woerdeman derived a semidefinite representable criterion for non-negativity of a partial finite Fourier series sum. Their criterion reads as follows:

\begin{tvr}\label{Lean-pos} \cite[Prop. 2.3]{mclean}
Let $\sigma_0\in\R, \sigma_k=\bar\sigma_{-k}\in\CC$ for $k=1, \cdots, N-1$. Then the finite Fourier series expansion  $\sigma(\nu) = \sum_{k=-N+1}^{N-1} \sigma_k e^{i k\nu}$ is a nonnegative function $\sigma(\nu)\ge 0$ for $\nu\in\R$, if and only if there exists a positive semidefinite Hermitian matrix  $F\in \mathcal{H}^N, F\succeq 0,$ and such that, for each $k=0, 1,\cdots, N-1$,
\[
\sum_{p=k+1}^N F_{p,p-k} = \sigma_k.
\]

\end{tvr}
Using Proposition~\ref{Lean-pos} and taking into account that $\sigma(\nu)+\sigma''(\nu) = \sum_{k=-N+1}^{N-1} (1-k^2) \sigma_k e^{i k \nu}$ for any $\sigma\in\K^N$ we end up with the following representation of the cone $\K^N$:

\begin{lem}\label{coneKrepresentation}
$\sigma\in\K^N$ if and only if there exist $F, G\in \mathcal{H}^N, F,G \succeq 0,$ such that, for any $k=0, \cdots, {N-1}$, 
\[
\sum_{p=k+1}^N F_{p,p-k} = \sigma_k, 
\quad
\sum_{p=k+1}^N G_{p,p-k} = (1-k^2) \sigma_k.
\]
\end{lem}

\section{Semidefinite programming (SDP) representability of the minimal anisotropic interface energy problem}

In this section we will show how the optimization problems \eqref{first} and \eqref{second} for construction of the optimal anisotropy function $\sigma$ can be reformulated in terms of the non-convex quadratic optimization problem \eqref{P}. In order to compute the function $\sigma$ we restrict ourselves to the finite dimensional subspace formed by anisotropy functions belonging to the finite dimensional cone $\K^N$ for given finite number of Fourier modes $N\in\N$. 

First we rewrite problems  \eqref{first} and \eqref{second} in terms of real and imaginary parts $x_\Re, x_\Im\in\R^N$ of the complex vector $\sigma\in\K^N$ representing the anisotropy function $\sigma$, i.~e. 
\[
x = \left(\begin{array}{c}
x_\Re \\ x_\Im
\end{array} \right)\equiv [x_\Re; x_\Im] \in \R^{n},\quad \sigma_k = x_{\Re,k} + i\,x_{\Im, k},
\]
for $k=0, \dots, N-1,$ where $n=2N$. Let $c_k, k=0, \dots, N-1$ be Fourier coefficients associated to the given Jordan curve $\Gamma$ (see \eqref{fouriercoeff}). If we set $\alpha= (\alpha_0, \cdots, \alpha_{N-1})^T, \beta=(\beta_0, \cdots, \beta_{N-1})^T \in \R^N$ where 
$\alpha_k = 2\Re c_k, \ \beta_k = 2 \Im c_k, k\ge1,  \ \alpha_0 = c_0, \beta_0=0$, then the anisotropic interface energy $L_\sigma(\Gamma)$ can be expressed as follows:
$L_\sigma(\Gamma) = \alpha^T x_\Re + \beta^T x_\Im$.

\subsection{SDP representation of problem \eqref{first} with linear constraints}

Using the representation of $L_\sigma(\Gamma)$ and semidefinite representation of $\K^N$ problem \eqref{first} with the linear constraint on $\sigma$ can be rewritten as optimization problem \eqref{P}

\begin{equation}\label{firstSDP}
\begin{array}{rl}
\displaystyle \min & 2 q_0^T x \\
{\rm s. t.} & \sigma_0 =1,
\\ & 
\sum_{p=k+1}^N F_{p,p-k} = \sigma_k,
\\ &
\sum_{p=k+1}^N G_{p,p-k} = (1-k^2)\sigma_k,
\\
& \hbox{for}\ \  k=0,\cdots, N-1,
\\ & F, G\succeq 0,
\end{array} 
\end{equation}
where $2 q_0 = [\alpha; \beta]\in \R^n, n=2 N$. It means that matrices $P_l=0$ for $l=0,1$ in \eqref{P}. In this case problem \eqref{P} is just a convex semidefinite programming problem with the linear value function and linear matrix inequality constraints. It can be solved directly by computational tools for solving SDP optimization problems over symmetric cones like e.~g. the Matlab software package SeDuMi by Sturm \cite{sturm}.

Notice that problem \eqref{firstSDP} is feasible because $x\equiv [x_\Re; x_\Im]$ where $x_{\Re,0}=1, x_{\Re,k}=0, k\ge 1, x_\Im =0$ is a feasible solution corresponding to the constant anisotropy function $\sigma\equiv 1$. Since $2q_0$ is the positive vector and $x\ge0$ for any $x$ feasible to  \eqref{firstSDP} we have $\hat p_1 >-\infty$.

\subsection{SDP representation of problem \eqref{second} with quadratic constraints}

According to the scaling property \eqref{homog} problem \eqref{second} with quadratic constraint imposed on $\sigma$ is equivalent (up to a positive scalar multiple of the optimal function $\sigma$) to the following finite dimensional optimization problem:
\begin{equation}\label{maxWN}
\begin{array}{rl}
\displaystyle \max_\sigma & |W_{\sigma}| \\
{\rm s. t.} & L_{\sigma}(\Gamma)=L(\Gamma),
\  \sigma \in \K^N, 
\end{array} 
\end{equation}
i.~e. we maximize the area $|W_\sigma|$ of the Wulff shape under the constraint that the interface energy $L_\sigma(\Gamma)$ is fixed to the predetermined constant, e.~g. the total length $L(\Gamma)$. The choice of the scaling constraint $L_\sigma(\Gamma)=L(\Gamma)$ is quite natural because in the case $\Gamma$ is a circle, the anisotropy function $\sigma\in \K$ maximizing $|W_\sigma|$ under the constraint  $L_\sigma(\Gamma)=L(\Gamma)$ is just unity, $\sigma\equiv1$.
Taking into account representation of the Wulff shape area from Section IV, part D, and introducing the real $n\times n$ matrix $P_0$:
\begin{equation}
P_0 =\hbox{diag}(p_0, p_1, \cdots, p_{N-1}, q_0, q_1, \cdots, q_{N-1}), 
\label{P0def}
\end{equation}
where $p_0=q_0=-\pi, p_k=q_k= 2\pi (k^2-1)$ for $k\ge 1$, the optimization problem \eqref{maxWN} can be rewritten as follows:
\begin{equation}\label{maxWNP}
\begin{array}{rl}
 \displaystyle \min & x^T P_0 x \\
{\rm s. t.} & A x = b, \qquad  
\\ & 
\sum_{p=k+1}^N F_{p,p-k} = \sigma_k,
\\ &
\sum_{p=k+1}^N G_{p,p-k} = (1-k^2)\sigma_k,
\\
& \hbox{for}\ \  k=0,\cdots, N-1,
\\ & F, G\succeq 0,\ \  F,G\in {\mathcal H}^N, x\in \R^{2N},
\end{array} 
\end{equation}
where $A$ is a $1\times n$ real matrix, $A= (\alpha^T, \beta^T)$ and $b=L(\Gamma)$. 
Since the matrix $P_0$ is indefinite problem then \eqref{maxWNP} is a non-convex optimization problem with LMI constraints having the form of \eqref{P}. Notice that the semidefinite constraints imposed on matrices $F$ and $G$ can be rewritten in terms of the LMI constraint in \eqref{P} by using a standard basis in the space of complex $N\times N$ Hermitian matrices. 

With regard to results from Section II the enhanced semidefinite relaxation of  problem \eqref{maxWNP} has the form:
\begin{equation}\label{maxWNP-relax}
\begin{array}{rl}
\displaystyle \min & \hbox{tr}( P_0 X) \\
{\rm s. t.} & A x = b,\  A X = b x^T, 
\\ & X \succeq x x^T, 
\\ & \sum_{p=k+1}^N F_{p,p-k} = \sigma_k, 
\\ & \sum_{p=k+1}^N G_{p,p-k} = (1-k^2)\sigma_k,
\\ & \hbox{for}\ k=0,\cdots, N-1, 
\\ & F, G\succeq 0,\  F,G\in {\mathcal H}^N, x\in \R^{2N}, X\in  {\mathcal S}^{2N}. 
\end{array} 
\end{equation}
Similarly as in the case of problem \eqref{firstSDP}, both problems \eqref{maxWNP} as well as \eqref{maxWNP-relax} are feasible because $x\equiv [x_\Re; x_\Im]$ where $x_{\Re,0}=1, x_{\Re,k}=0, k\ge 1, x_\Im =0$ is a feasible solution to  \eqref{maxWNP} and $(x, X), X=xx^T,$ is feasible to  \eqref{maxWNP-relax}. Moreover, the optimal value $\hat p_1$ is finite. Indeed, it follows from the anisoperimetric inequality $\Pi_\sigma(\Gamma)\ge 1$ that 
$- x^T P_0 x = |W_\sigma| = \frac{L_\sigma(\Gamma)^2}{4\Pi_\sigma(\Gamma){\mathcal A}(\Gamma)}\le \frac{L(\Gamma)^2}{4 {\mathcal A}(\Gamma)} <\infty,$
for any $\sigma = x_\Re + i\,x_\Im\in\K^N,\ x=[x_\Re; x_\Im],$ which is a feasible solution to \eqref{maxWNP}. Hence $\hat p_1> -\infty$. 

\section{Numerical experiments based on minimization of the interface energy}

\subsection{The case of a linear constraint}

First, we present results of numerical resolution of the optimal anisotropy function $\sigma$ based on a solution to problem \eqref{first}. It should be obvious that, up to a positive multiple of $\sigma$, the optimal solution $\sigma$ to the maximization problem \eqref{maxWN} is also an optimal solution to the minimal interface energy problem with the linear constraint $\sigma_{avg}\equiv \frac{1}{2\pi}\int_0^{2\pi} \sigma(\nu)\dnu  =\sigma_0 = 1$ imposed on $\sigma\in \K^N$, i.~e.
\begin{equation}\label{minLavg}
\begin{array}{rl}
\displaystyle\min_\sigma & L_\sigma(\Gamma) \\
{\rm s. t.} & \sigma_0 = 1,\ \  \sigma \in \K^N .
\end{array} 
\end{equation}
Notice that a minimizer $\sigma$ to \eqref{minLavg} need not be unique. Indeed, let $\Gamma$ be a circle with a radius $r>0$. Then the tangent angle $\nu$ can be used for parameterization of a convex curve $\Gamma$  such that  $\dnu = \kappa \ds = r^{-1} \ds$ (cf. \cite{MS2004,SY}). We obtain
\[
L_\sigma(\Gamma) = \int_\Gamma \sigma(\nu) \ds = r \int_0^{2\pi} \sigma(\nu) \dnu = 2\pi r,  \sigma_{avg} = 2\pi r
\]
for any $\sigma\in\K^N$ such that $\sigma_{avg} = 1$. But this means that any $\sigma\in\K$ with  $\sigma_{avg} = 1$ is a minimizer to \eqref{minLavg}.

\begin{figure}
\begin{center}
\subfigure[]{\includegraphics[width=0.24\textwidth]{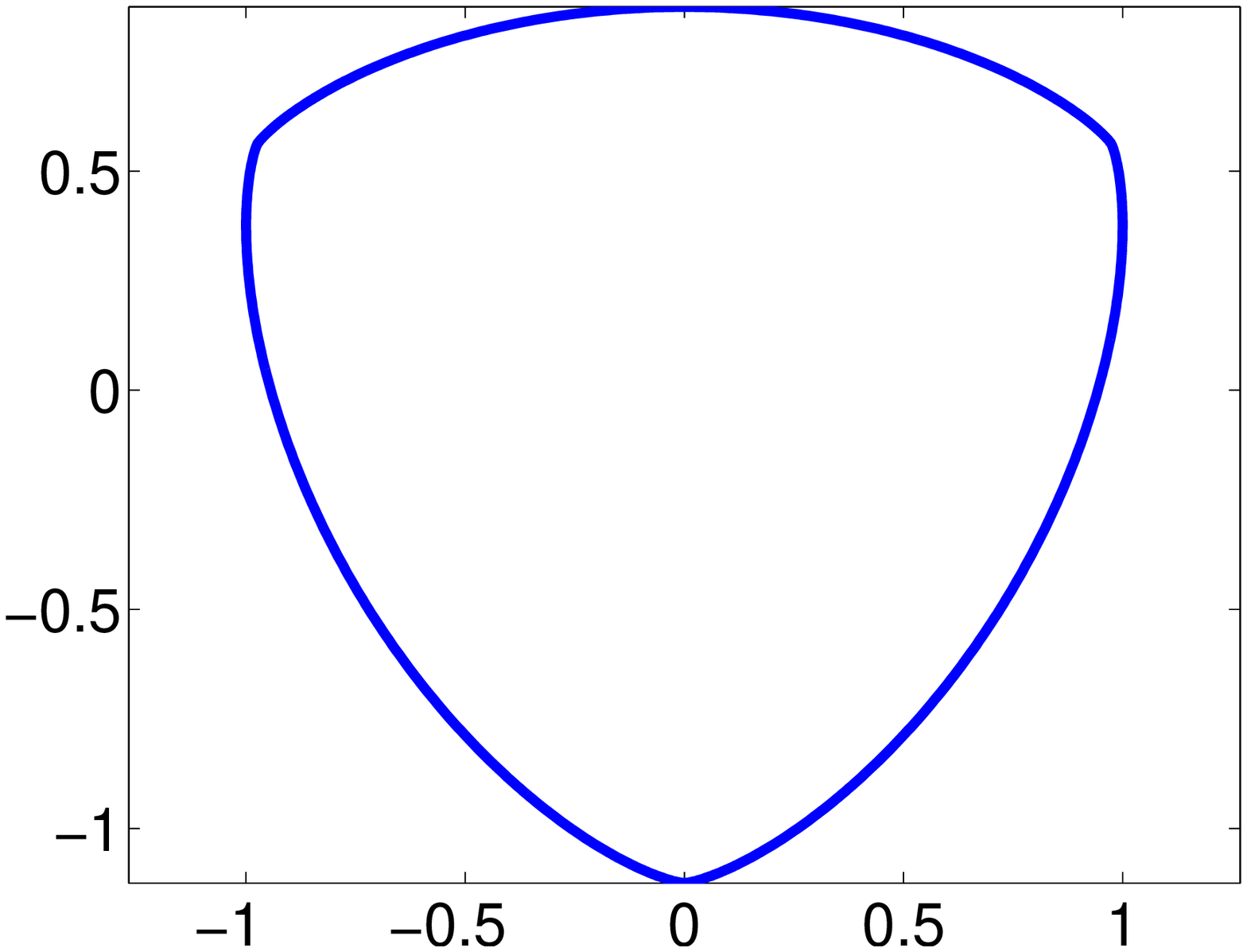}}
\subfigure[]{\includegraphics[width=0.24\textwidth]{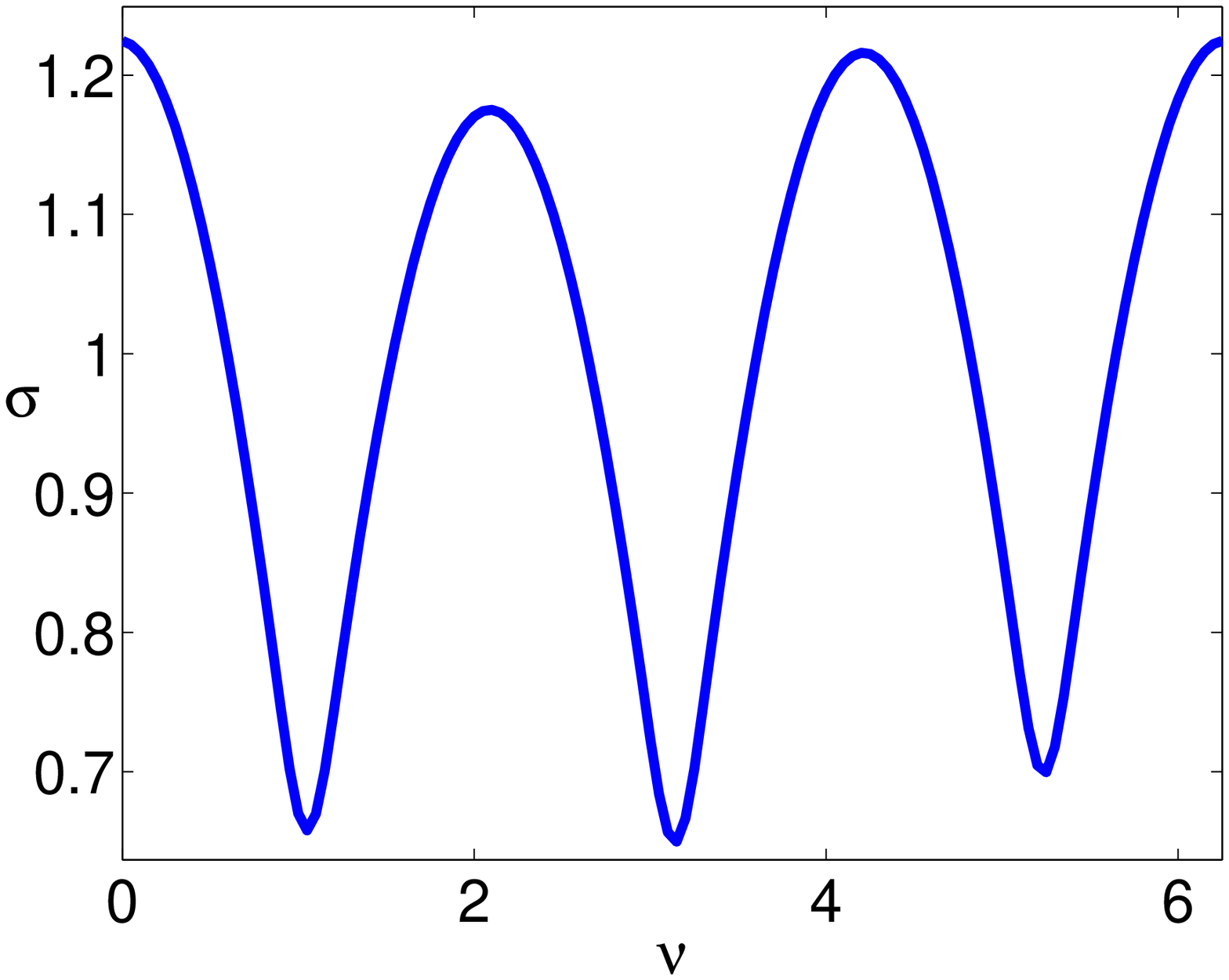}}
\\
\subfigure[]{\includegraphics[width=0.24\textwidth]{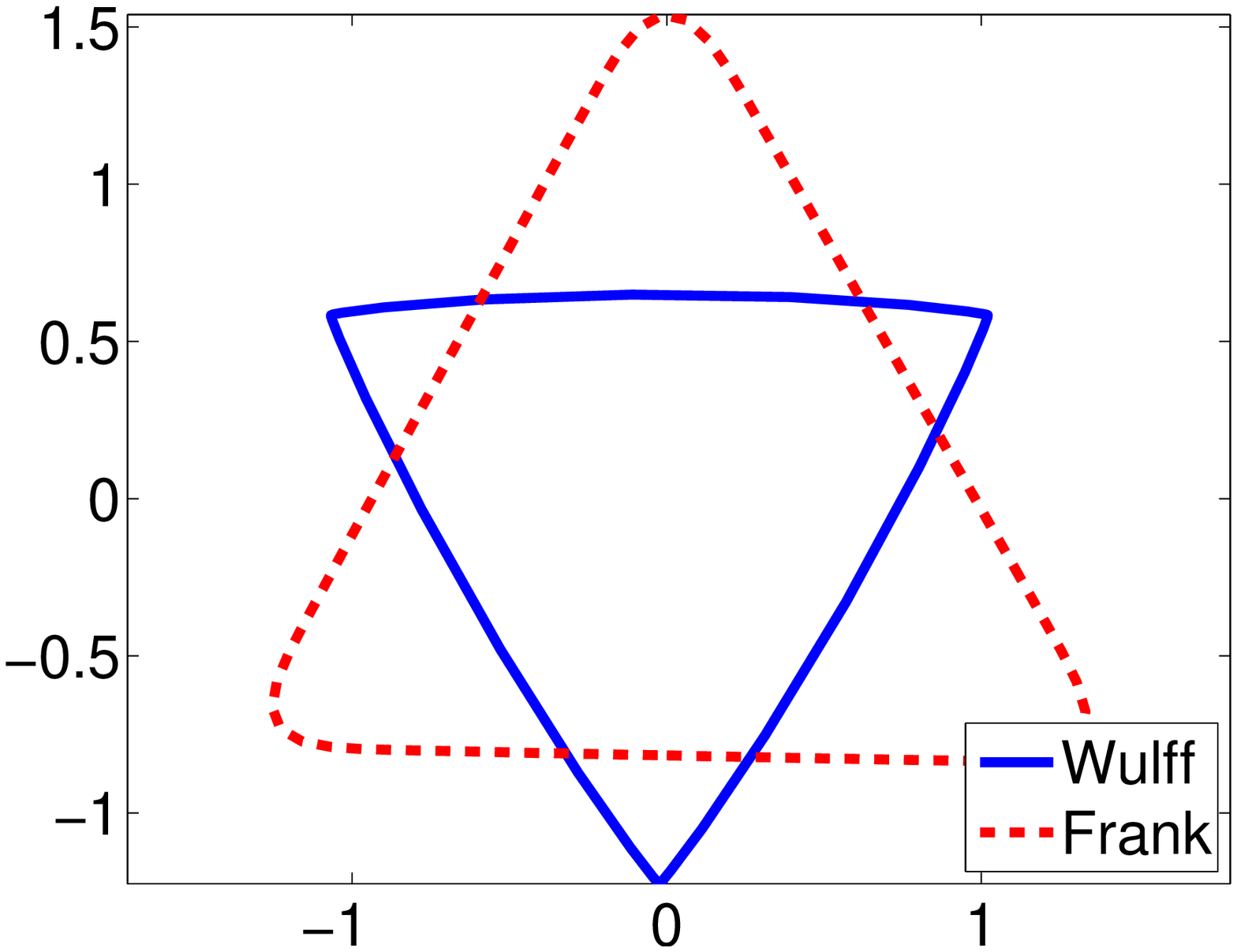}}
\subfigure[]{\includegraphics[width=0.24\textwidth]{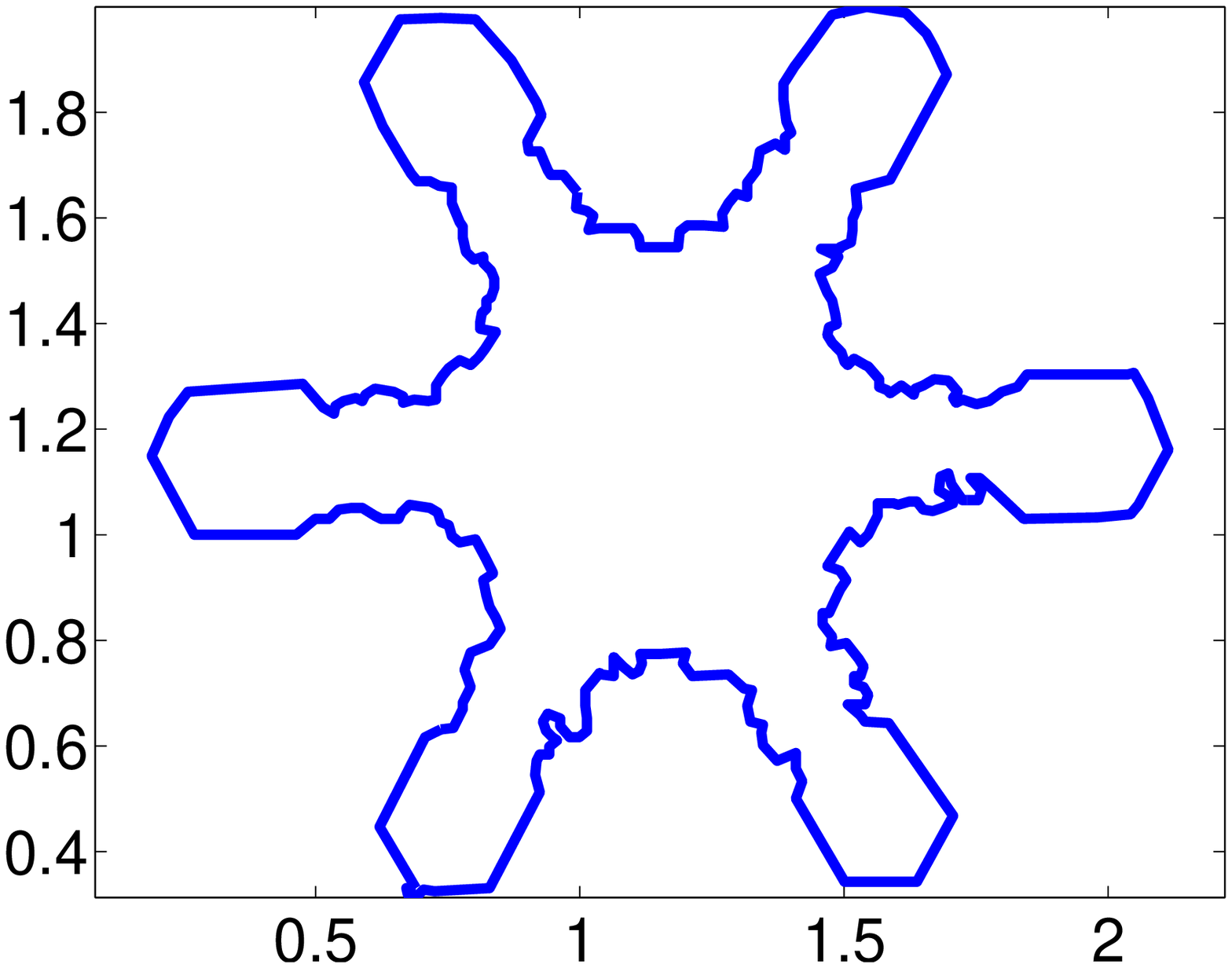}}
\\
\subfigure[]{\includegraphics[width=0.24\textwidth]{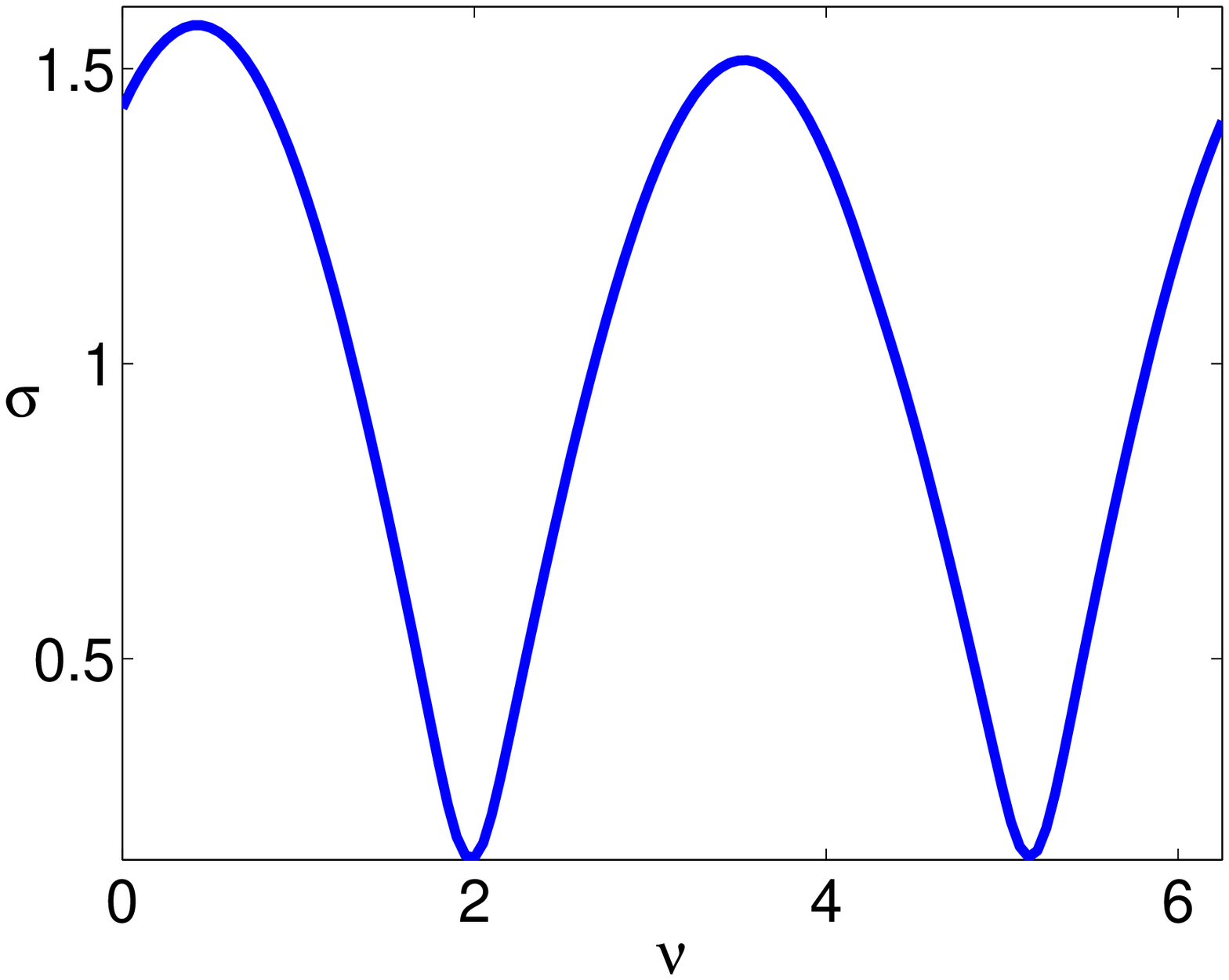}}
\subfigure[]{\includegraphics[width=0.24\textwidth]{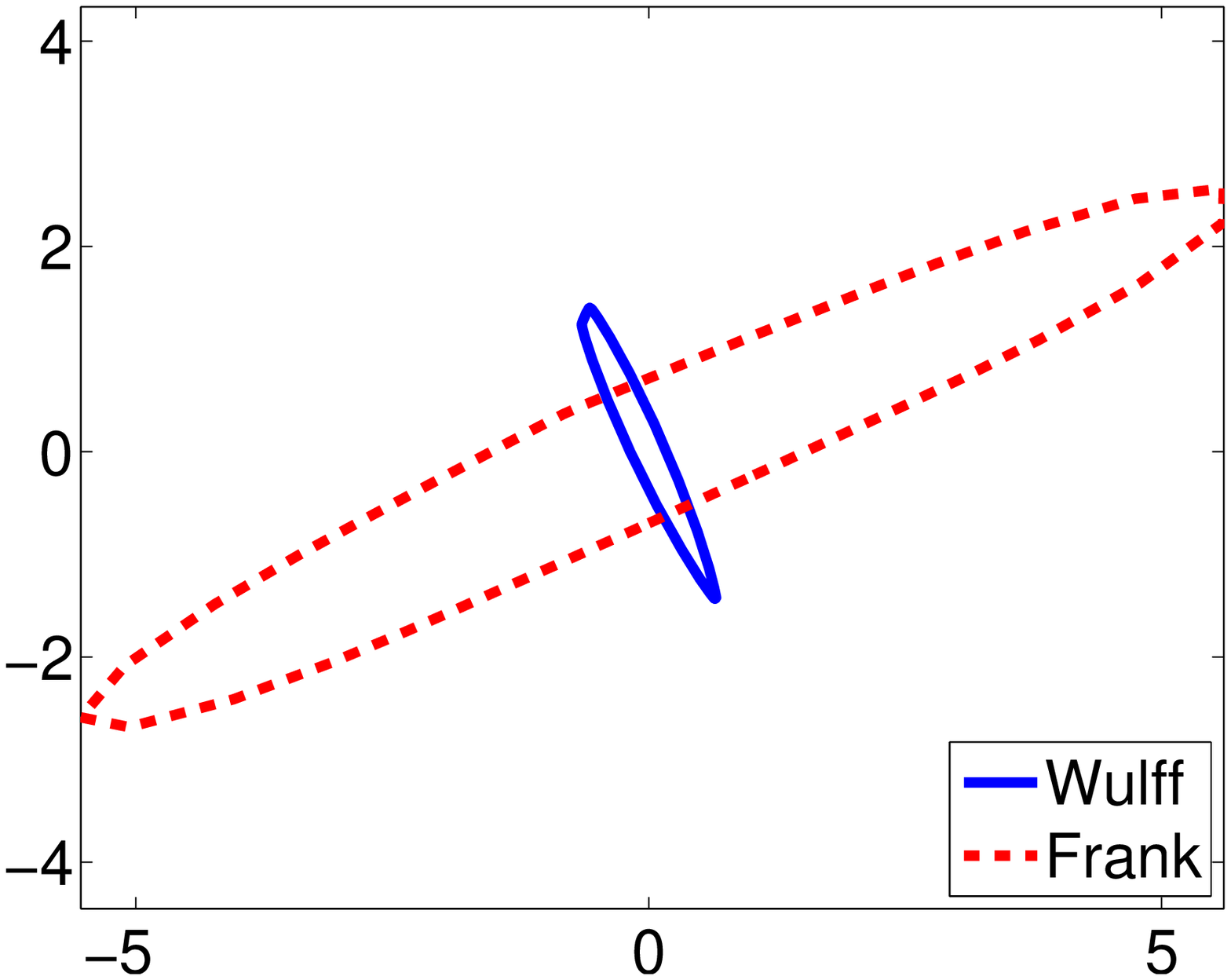}}

\end{center}
\caption{
Jordan curves corresponding to boundaries of the Wulff shape $W_\mu$ (a) and a snowflake (d). The optimal anisotropy functions $\sigma\in \K^N$ computed by means of a solution to the SDP  \eqref{first} are shown in (b) and (e), respectively. The Wulff shape and Frank diagram are depicted in (c) and (f)
}
\label{fig:wulff}
\end{figure}

In Fig.~\ref{fig:wulff} we plot two examples of Jordan curves $\Gamma$ in the plane. In  Fig.~\ref{fig:wulff} (a) we plot a Jordan curve representing the boundary of the Wulff shape $\partial_\mu W$. In this example the function $\mu$ is the Kobayashi three-fold anisotropy function $\mu(\nu) = 1+\varepsilon \cos(m\nu)$ with $m=3$ and $\varepsilon =0.99/(m^2-1)$. Clearly, $\mu\in \K^N$ for any $N\ge4$. Unfortunately, resolution of the optimal $\sigma$ based on a solution to \eqref{minLavg}, i.~e. \eqref{first} does not recover the original anisotropy function $\mu$ as one may expect in this case. The optimal anisotropy function $\sigma$ has the flattened Wulff shape having sharp corners. We also plotted the corresponding Frank diagram defined as follows:
\[
{\mathcal F}_\sigma = \{\vecx= - r \vecn \ |\  0\le r\le 1/\sigma(\nu), \nu\in [0,2\pi] \},
\]
where $\vecn=(-\sin\nu, \cos\nu)^T$. The next example of a curve $\Gamma$ representing a boundary of a real snowflake shown in Fig.~\ref{fig:wulff} (d) is even worse. The optimal anisotropy function $\sigma$ obtained by solving  \eqref{minLavg} has two local maxima corresponding thus to the two-fold anisotropy rather than hexagonal one, as one may expect in this case.

\subsection{The case of a quadratic constraint}

In this section we present results of resolution of the optimal anisotropy function by means of a solution to \eqref{second} in which $\sigma$ is minimizer of the anisotropic energy $L_\sigma(\Gamma)$ subject to the quadratic constraint $|W_\sigma|=1$. As it was already discussed in Section V, the optimization problem \eqref{second}  leads to a non-convex SDP problem \eqref{maxWN} which we can solve by means of the enhanced semidefinite relaxation problem \eqref{maxWNP-relax}. It was solved numerically by using the powerful nonlinear convex programming Matlab solver SeDuMi developed by J.~Sturm \cite{sturm}. Notice that it implements self-dual embedding method proposed by Ye, Todd and Mizuno \cite{todd}. It is worth to note that without complementing \eqref{maxWN} by the quadratic-linear constraint $AX=b x^T$ in \eqref{maxWNP-relax} the SeDuMi solver was unable to solve the problem because of its unboundedness. 
In order to call SeDuMi solver we have uased the CVX Matlab programming framework (cf. Henrion et al. \cite{henrion})




\begin{figure}
\begin{center}
\subfigure[]{\includegraphics[width=0.24\textwidth]{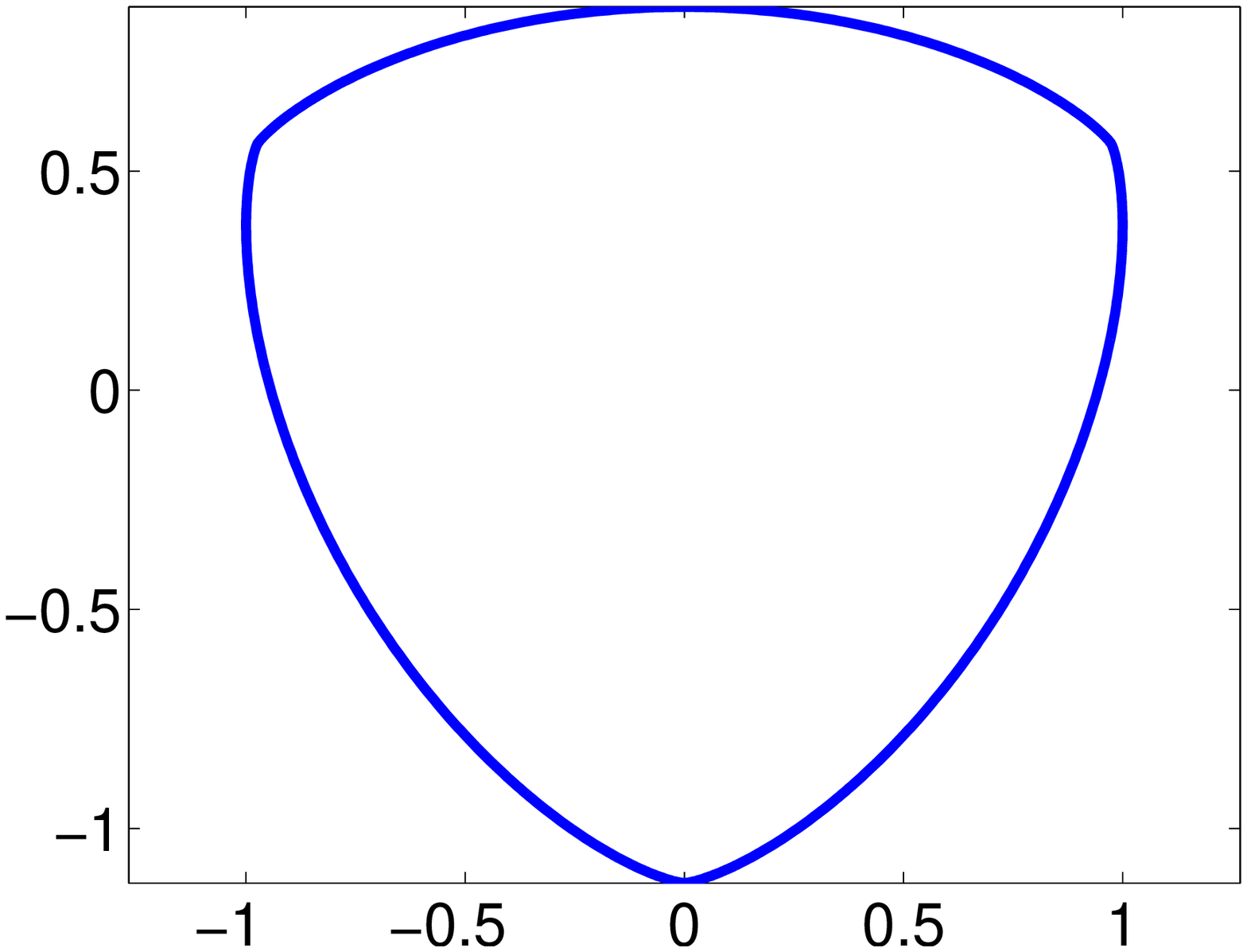}}
\subfigure[]{\includegraphics[width=0.24\textwidth]{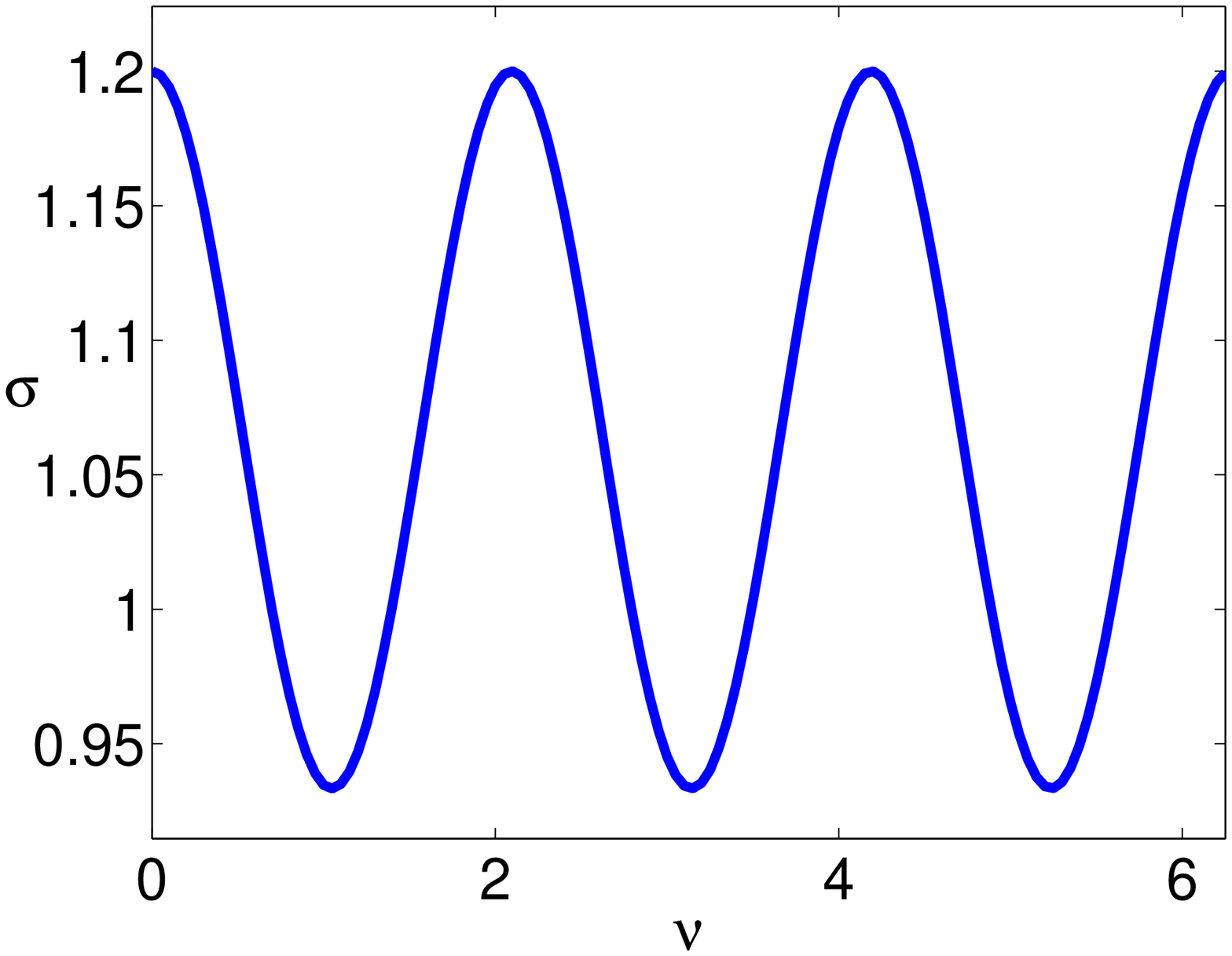}}
\\
\subfigure[]{\includegraphics[width=0.24\textwidth]{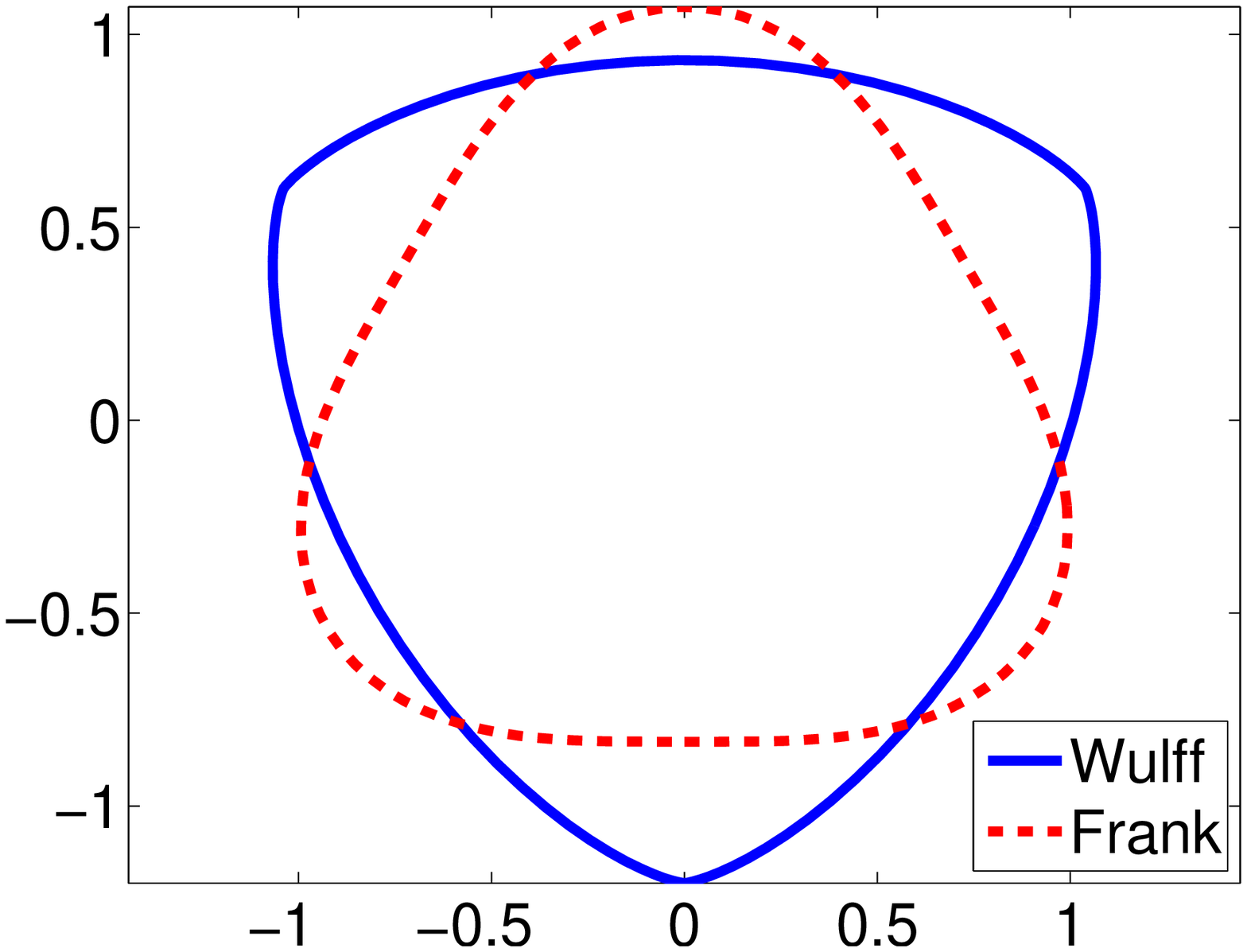}}
\end{center}
\caption{
A Jordan curve corresponding to the boundary of the Wulff shape $W_\mu$ (a). The optimal anisotropy functions $\sigma\in \K^N$ computed by means of a solution to the enhanced semidefinite relaxed program \eqref{maxWNP-relax} are shown in (b). The Wulff shape and Frank diagram are depicted in (c)
}
\label{fig:wulff2}
\end{figure}

\begin{figure}
\begin{center}
\subfigure[]{\includegraphics[width=0.24\textwidth]{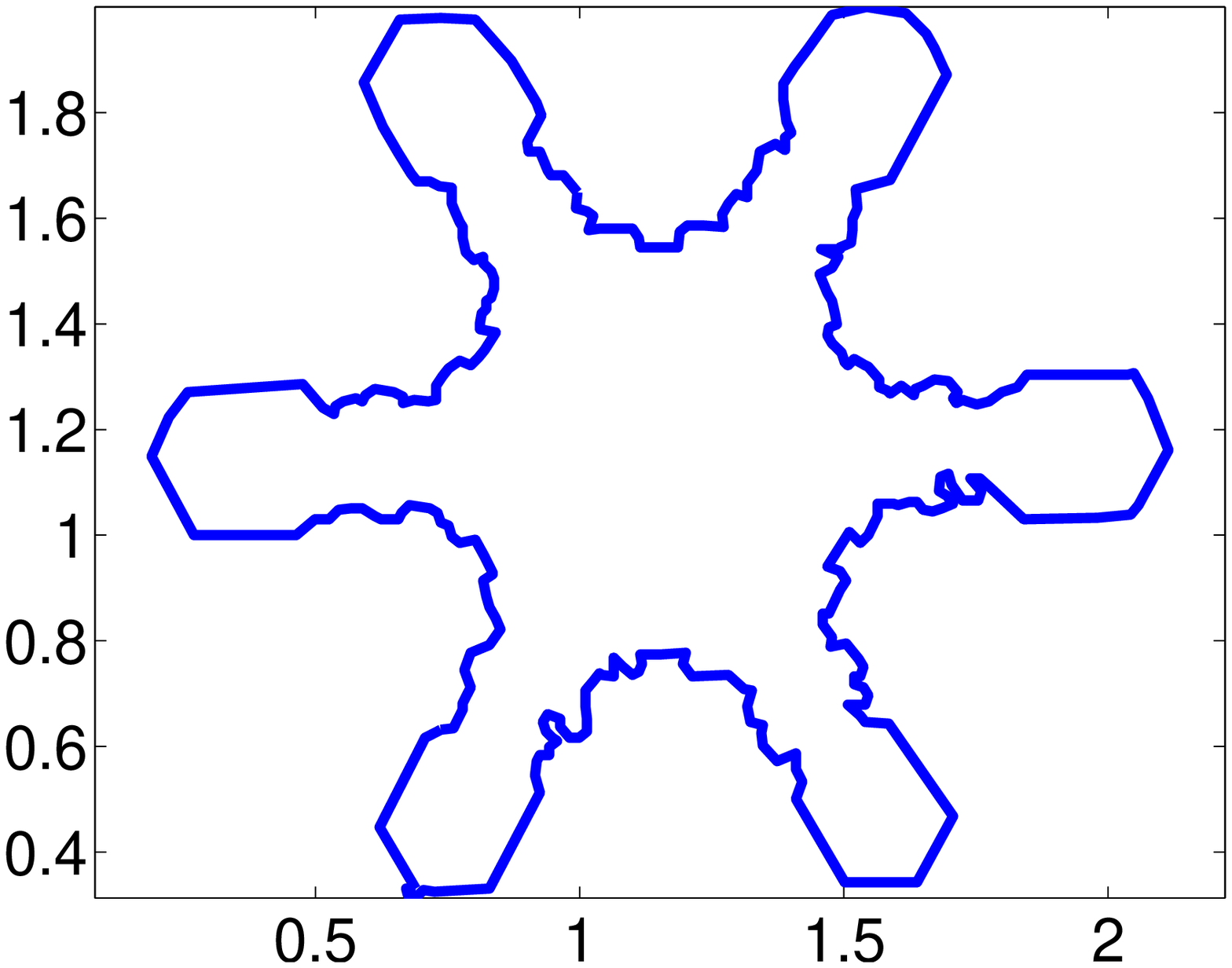}}
\subfigure[]{\includegraphics[width=0.24\textwidth]{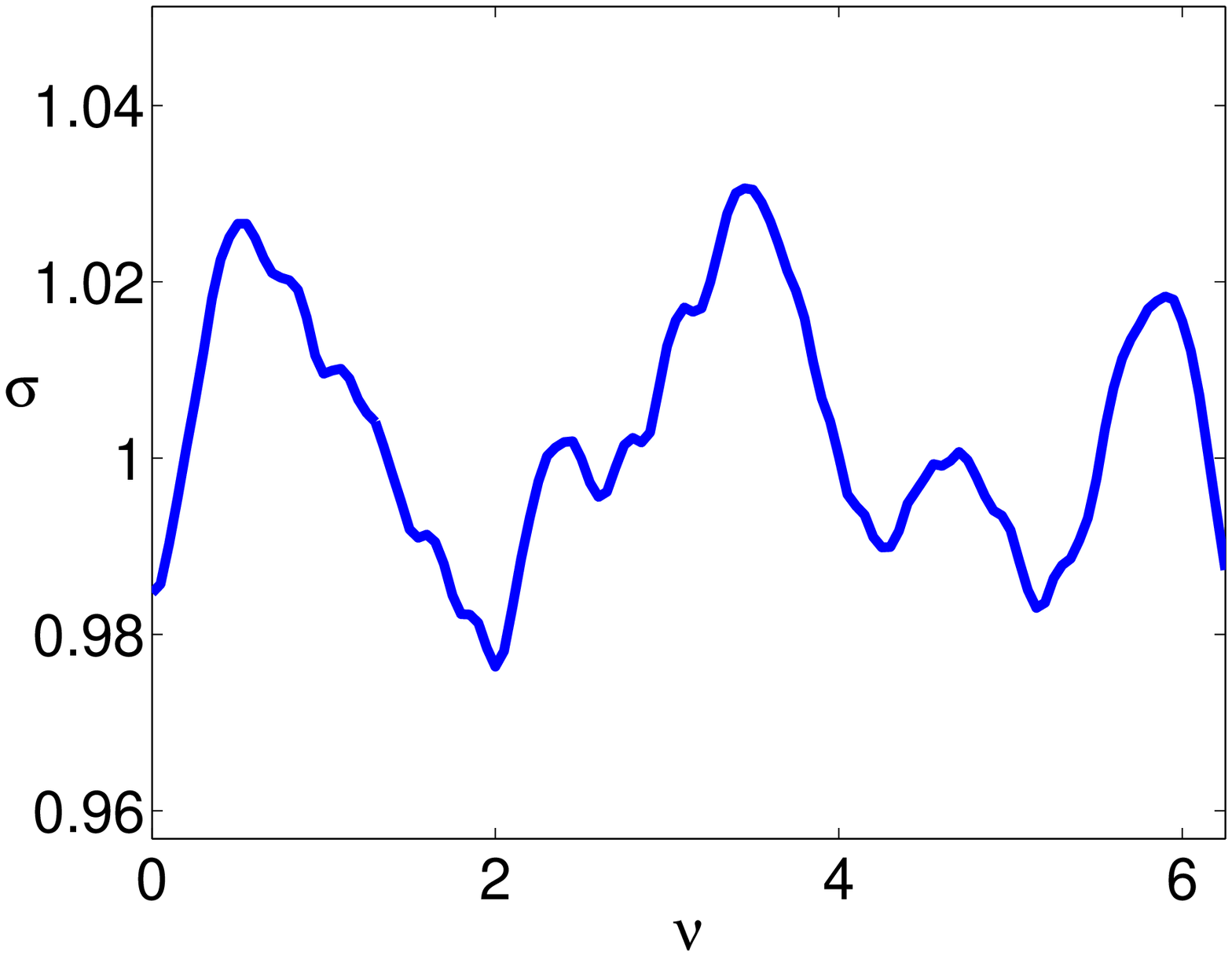}}
\\
\subfigure[]{\includegraphics[width=0.24\textwidth]{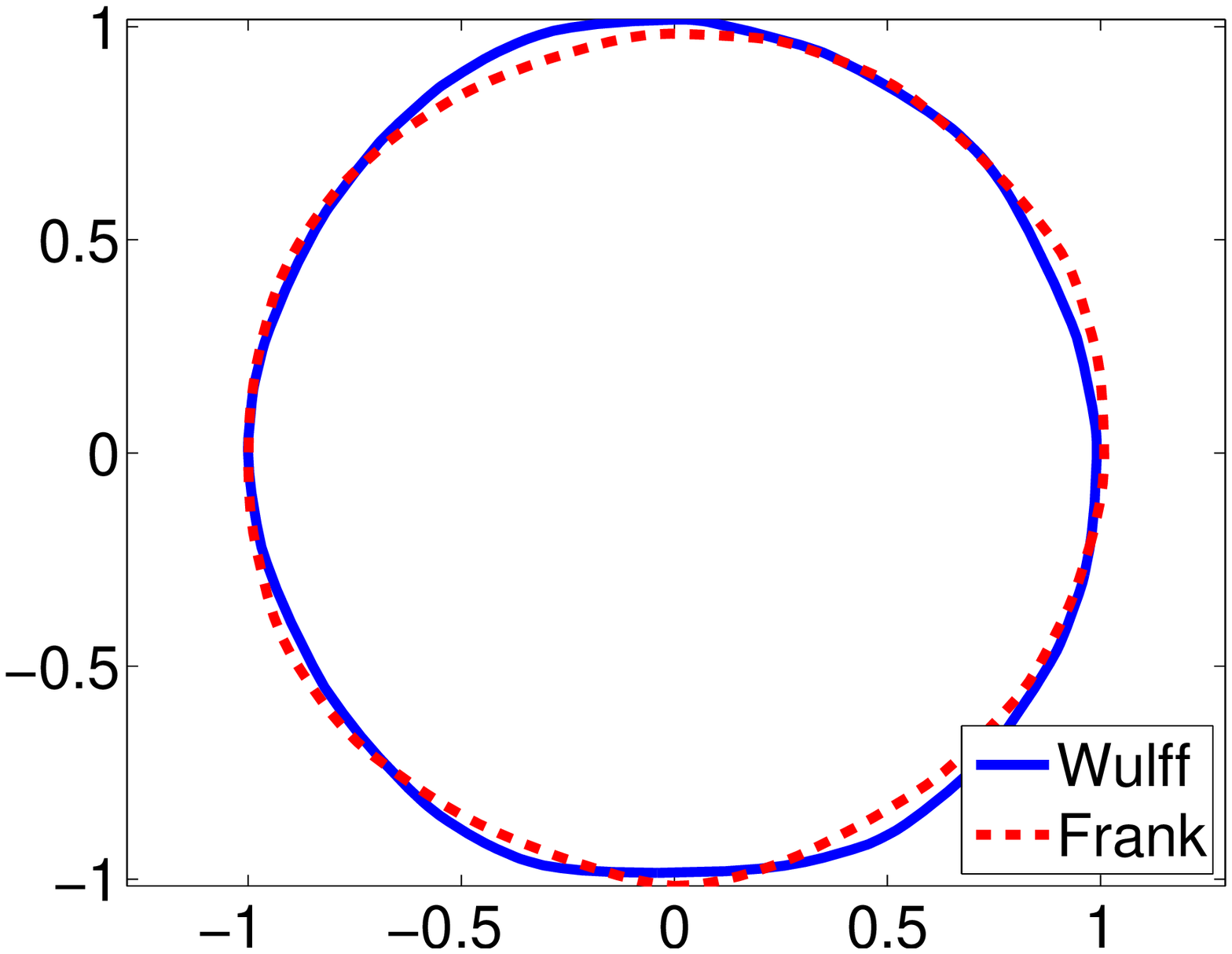}}
\subfigure[]{\includegraphics[width=0.24\textwidth]{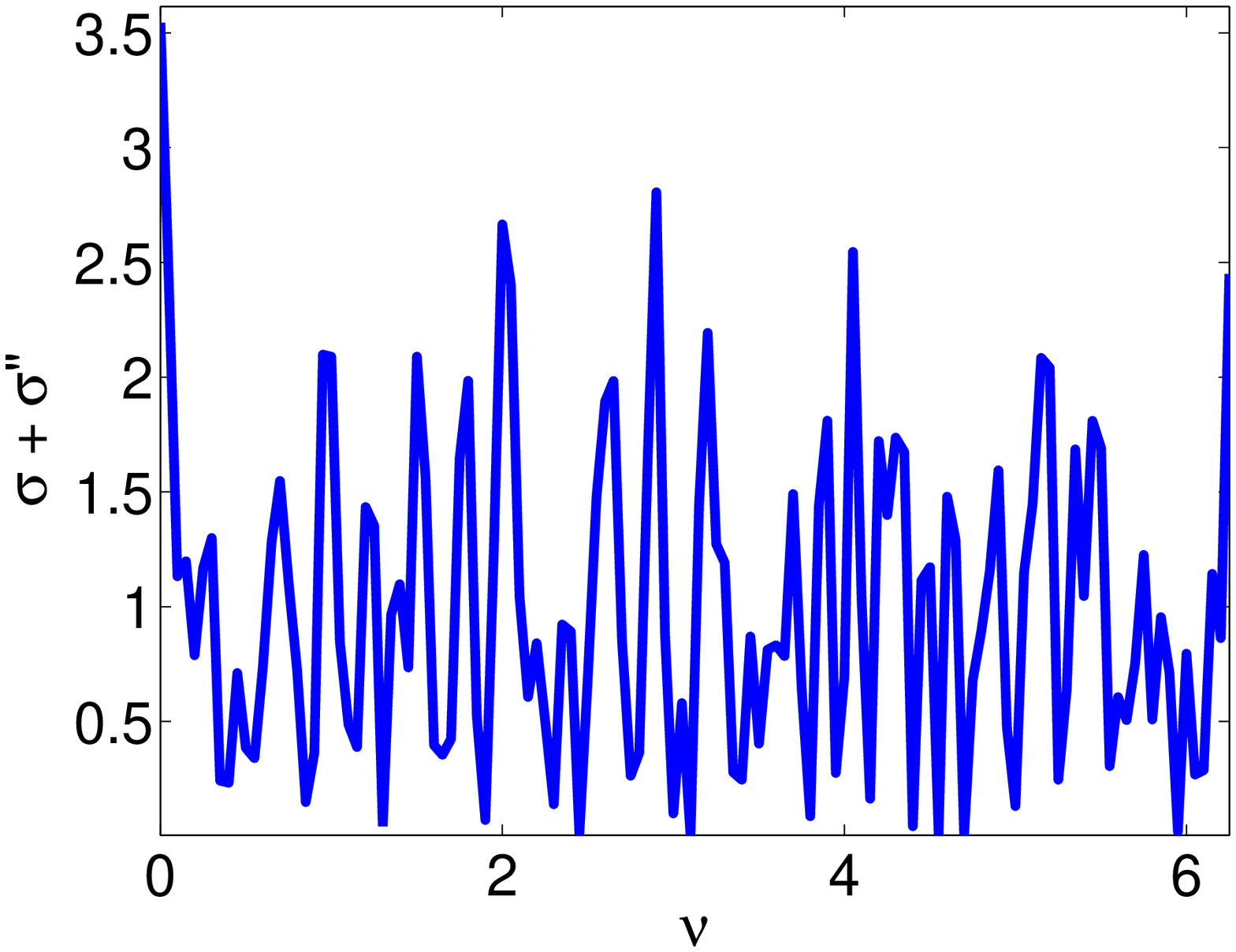}}
\end{center}
\caption{
A Jordan curve corresponding to the boundary of a snowflake (a). The optimal anisotropy function $\sigma\in \K^N$ computed by means of a solution to \eqref{maxWNP-relax} is shown in (b). The Wulff shape and Frank diagram are depicted in (c). The reciprocal value of the curvature $\kappa^{-1}= \sigma + \sigma''$ is depicted in (d)
}
\label{fig:snowflake4}
\end{figure}

In Fig.~\ref{fig:wulff2} we plot the same three-fold  Jordan curve $\Gamma$ as in Fig.~\ref{fig:wulff} (a). Using the quadratic constraint on $\sigma$ (see \eqref{second} and \eqref{maxWNP-relax}) the optimal anisotropy function $\sigma$ coincides with the Kobayashi three-fold anisotropy function $\mu(\nu) = 1+\varepsilon \cos(m\nu)$ (see (a,b,c)) with $m=3$ (cf. \cite{kobayashi}). In the case of a real snowflake boundary shown in Fig.~\ref{fig:snowflake4} (a), resolution of the optimal anisotropy function yields the Wulff shape as it can be seen from Fig.~\ref{fig:snowflake4} (c). There is just a small relative deviation (less than 2\%) of the function $\sigma(\nu)$ from the constant value $\sigma\equiv 1$. The same phenomena is however true for the Kobayashi function  $\mu(\nu) = 1+\varepsilon \cos(m\nu)$ with $m=6$ and $\varepsilon =0.99/(m^2-1)$.
The behavior of the optimal anisotropy can be better observed from Fig.~\ref{fig:snowflake4} (d) in which we plot the reciprocal value of the curvature $\kappa^{-1}= \sigma(\nu) + \sigma''(\nu)$ of the optimal Wulff shape. It has at least six separated spots of local minima close to zero value corresponding to high values of the curvature $\kappa$.

\begin{figure}
\begin{center}
\subfigure[]{\includegraphics[width=0.24\textwidth]{figures/dendride-curve.eps}}
\subfigure[]{\includegraphics[width=0.24\textwidth]{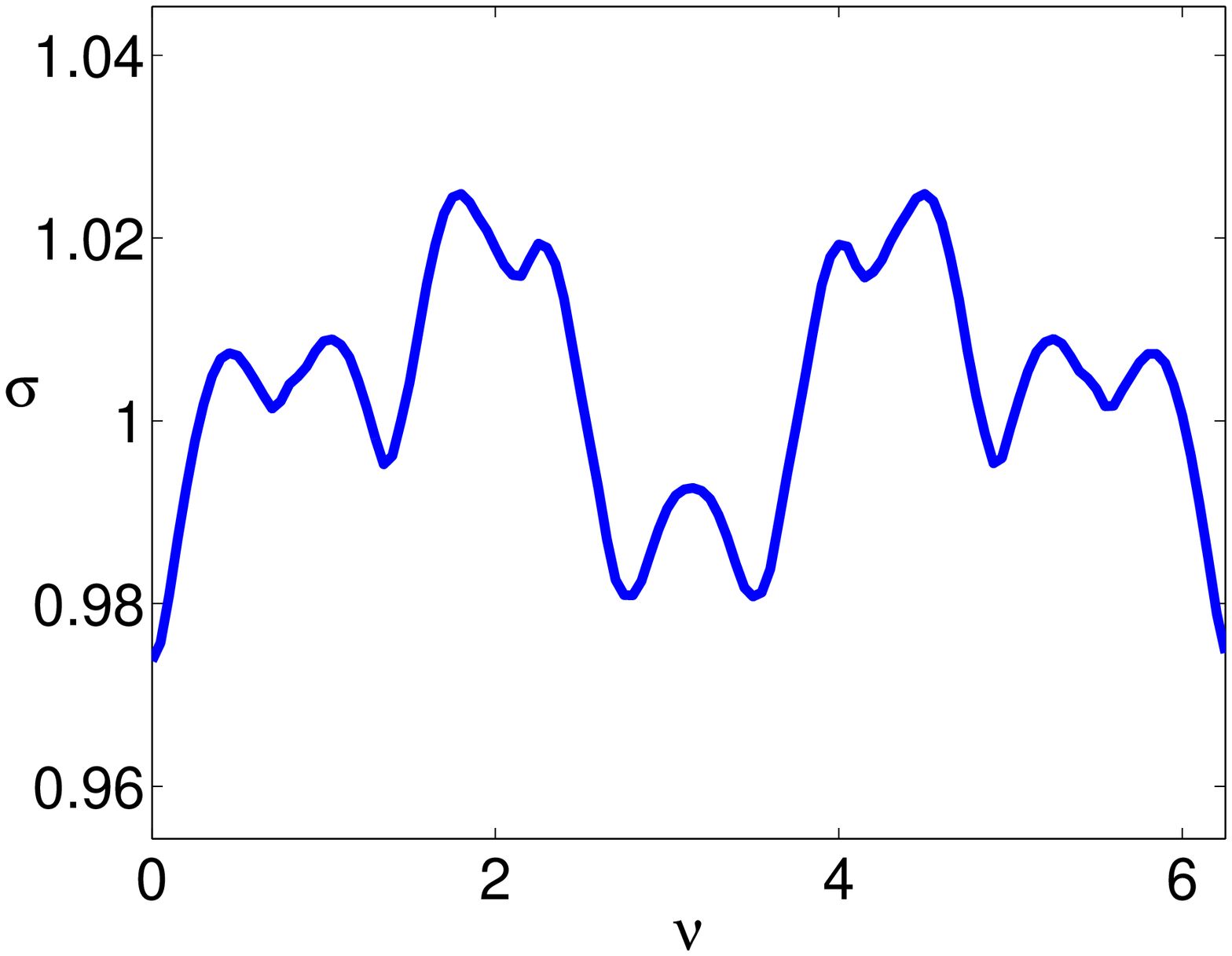}}
\\
\subfigure[]{\includegraphics[width=0.24\textwidth]{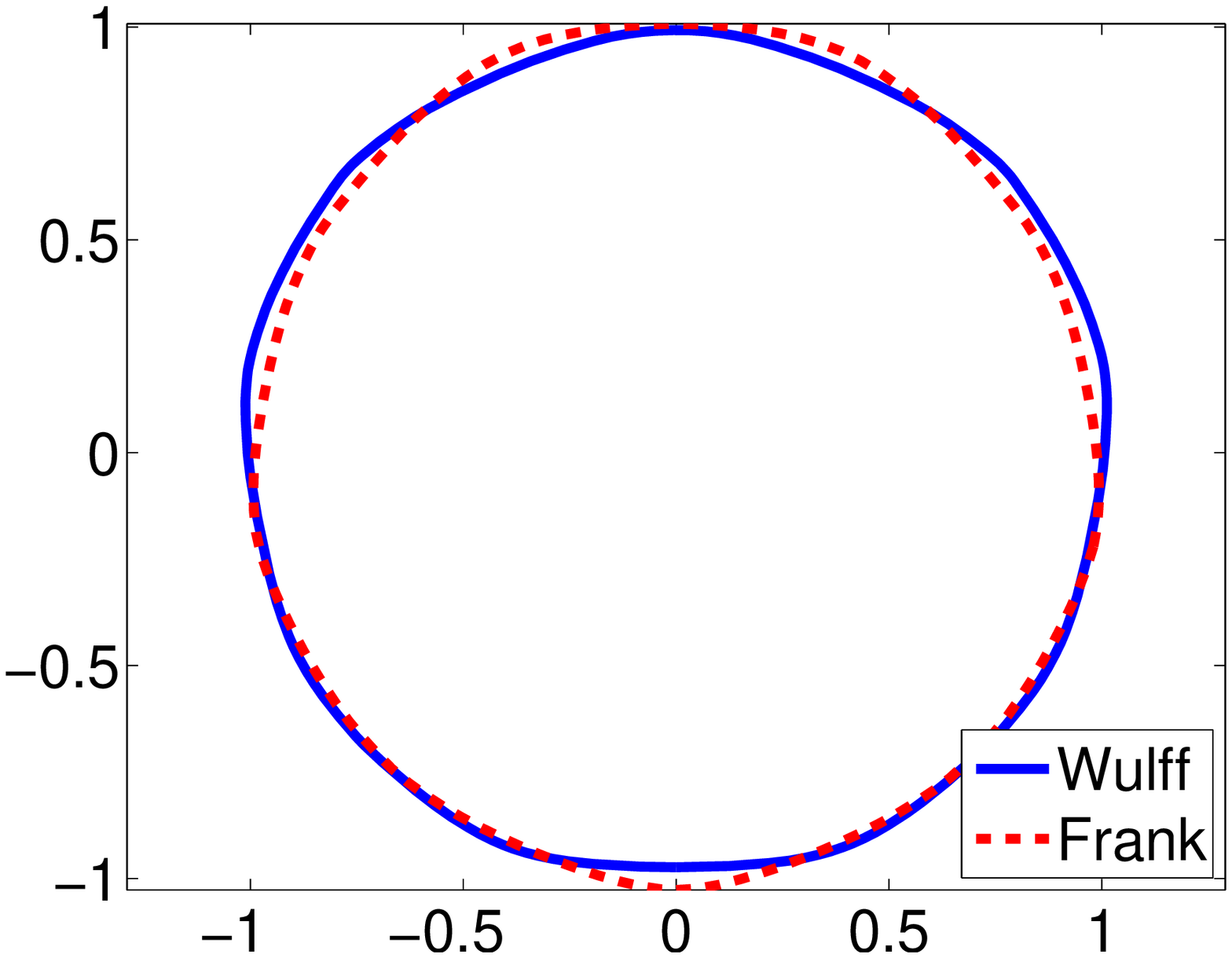}}
\subfigure[]{\includegraphics[width=0.24\textwidth]{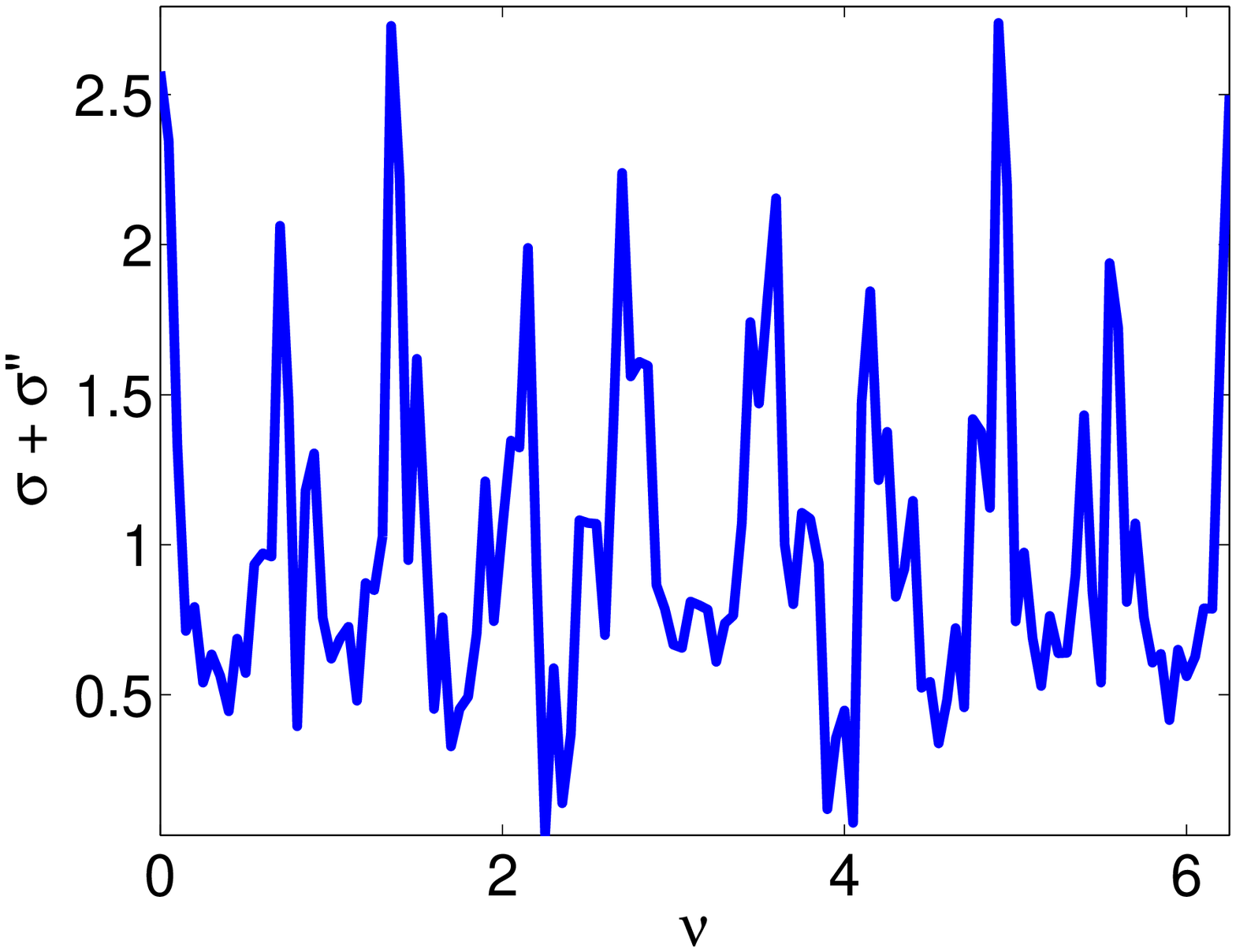}}
\end{center}
\caption{\small
A curve $\Gamma$ (a); the optimal anisotropy function $\sigma\in \K^N, N=50$ (b); the Wulff shape $W_\sigma$ and Frank diagram ${\mathcal F}_\sigma$ (c); the reciprocal value of the curvature $\kappa^{-1}= \sigma + \sigma''$ (d)
}
\label{fig:dendride}
\end{figure}

\begin{figure}
\begin{center}
\subfigure[]{\includegraphics[width=0.24\textwidth]{figures/snowflake6-curve.eps}}
\subfigure[]{\includegraphics[width=0.24\textwidth]{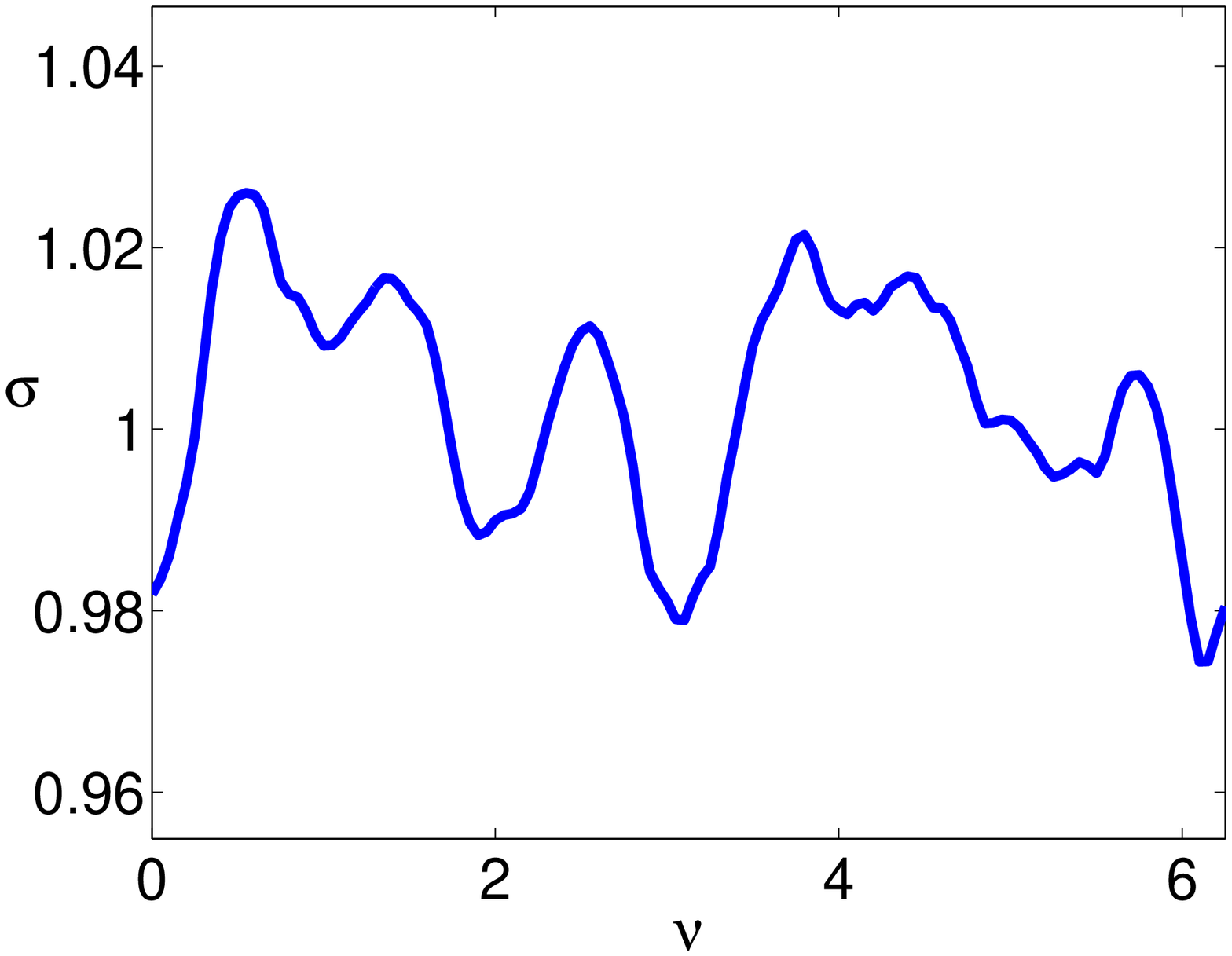}}
\\
\subfigure[]{\includegraphics[width=0.24\textwidth]{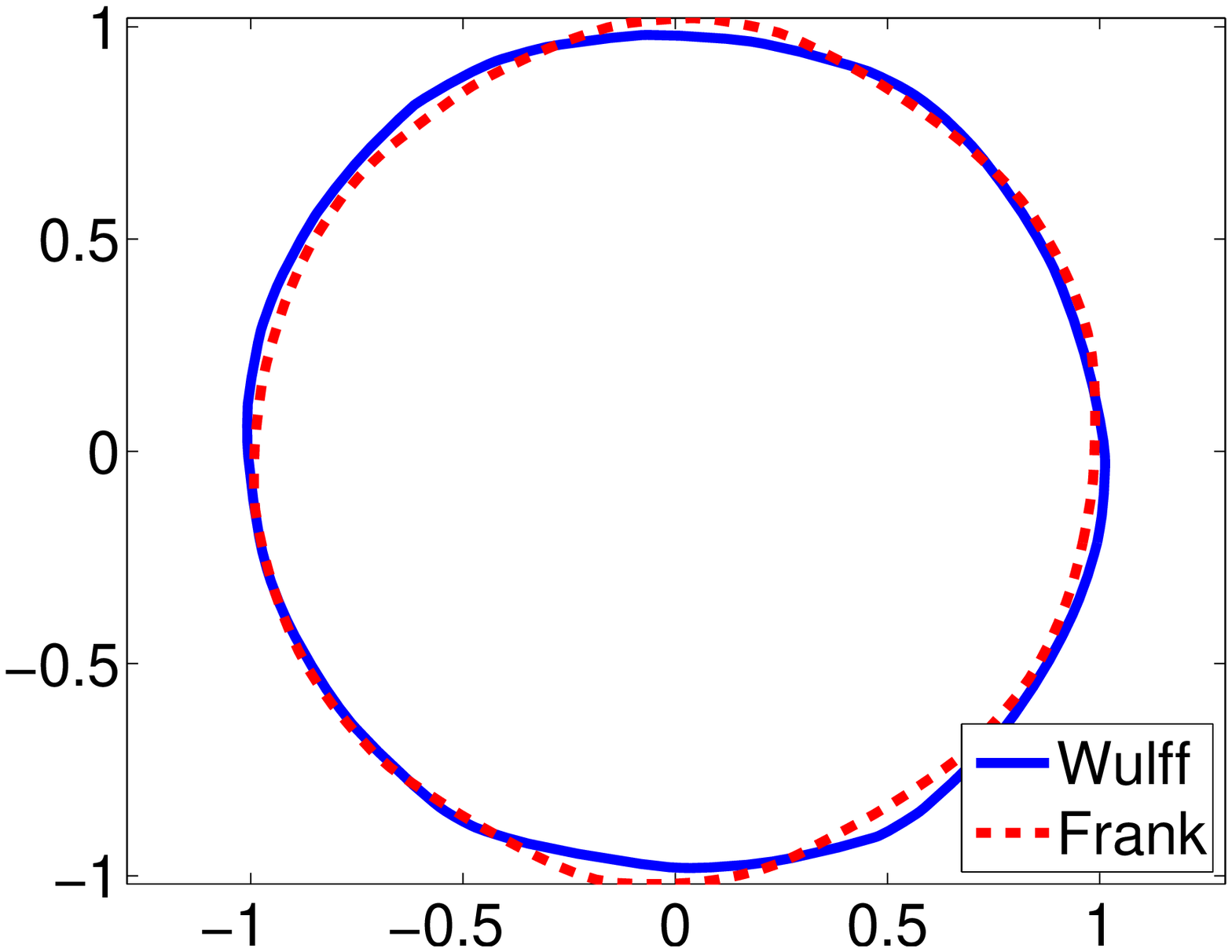}}
\subfigure[]{\includegraphics[width=0.24\textwidth]{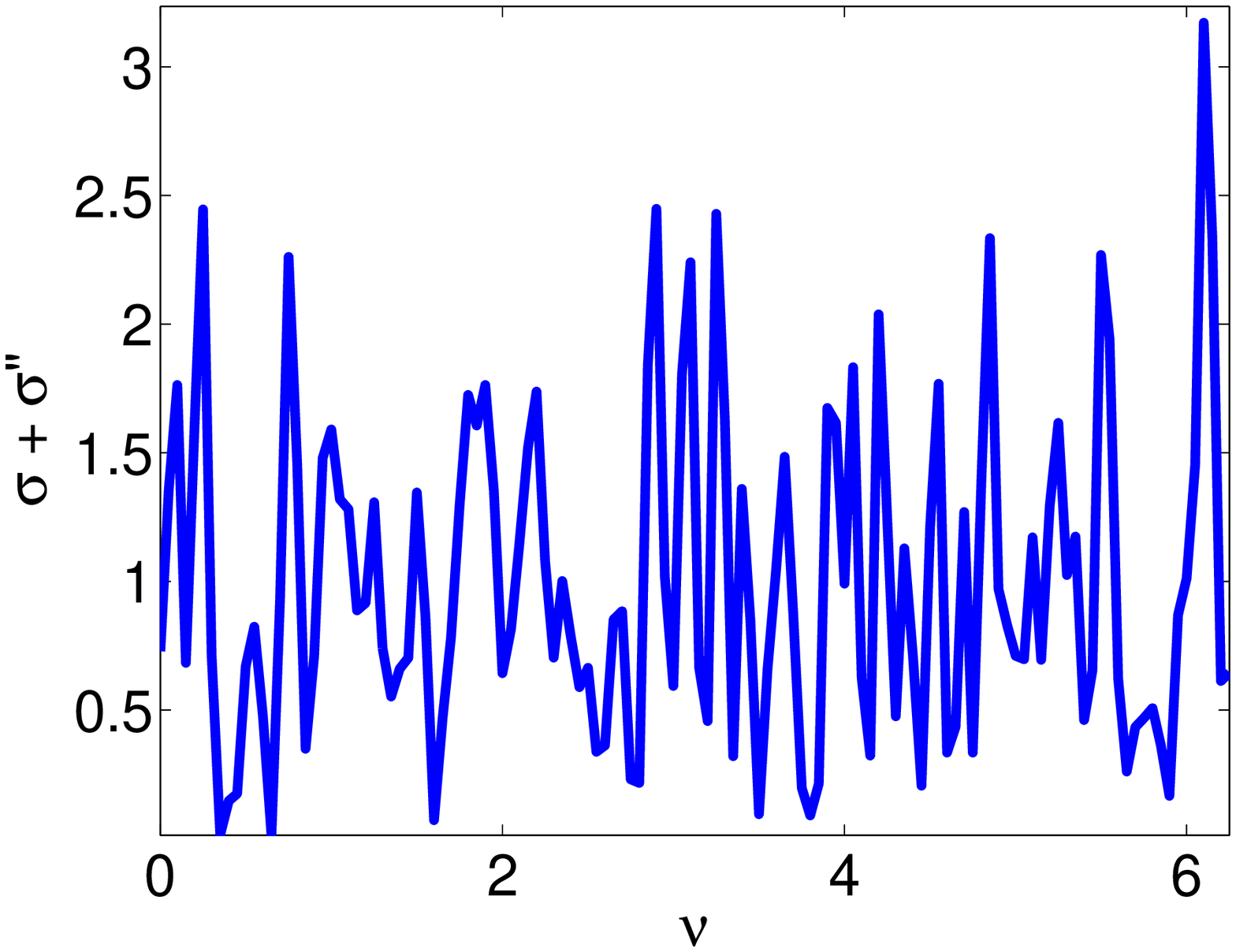}}
\end{center}
\caption{\small
A curve $\Gamma$ (a); the optimal anisotropy function $\sigma\in \K^N, N=50$ (b); the Wulff shape $W_\sigma$ and Frank diagram ${\mathcal F}_\sigma$ (c); the reciprocal value of the curvature $\kappa^{-1}= \sigma + \sigma''$ (d)
}
\label{fig:snowflake6}
\end{figure}

In Fig.~\ref{fig:dendride} (a) we present a simple test example of a Jordan curve $\Gamma$ given by the parameterization:  $\Gamma=\{\vecx(u)=r(u)(\sin(2\pi u), \cos(2\pi u))^T \ | \ u\in[0,1]\}$ where $r(u)=3+\exp(\cos(18\pi u))\cos(8\pi u)$.
The curve $\Gamma$ has been discretized by $K=1000$ grid points and the Fourier coefficients were computed according to \eqref{cp-discr}. We chose $N=50$ Fourier modes in this example. 
In the last numerical example shown in Fig.~\ref{fig:snowflake6} (a) we present computation of the optimal anisotropy function $\sigma$ for a boundary $\Gamma$ of a real snowflake. We again used $N=50$ Fourier modes and $K=700$ grid points for approximation of the boundary of a  snowflake. The resulting optimal anisotropy function $\sigma$ again corresponds to the Wulff shape with hexagonal symmetry. It can be seen from the plot of Fig.~\ref{fig:snowflake6} (d) in which we can observe six distinguished local minima of the reciprocal value $\kappa^{-1}$ of the curvature. 

In Table~\ref{tab1} we present results of numerical computations for various numbers $N$ of Fourier modes for the curve $\Gamma$ shown in Fig~\ref{fig:dendride} (a). We calculated the experimental order of time complexity (eotc) by comparing elapsed times $T_k$ for different $N_k$ as follows: $eotc_k = \ln(T_{k+1}/T_k)/\ln(N_{k+1}/N_k)$. It turns out that the time complexity measured by the $eotc$ is  below the order of $4.5$. On the other hand, for practical purposes, taking $N\approx 100$ Fourier modes is sufficient. Numerical computations were performed on a Quad-Core AMD Opteron Processor with 2.4GHz frequency, 32GB of memory. We also computed the relative gap in the optimal solution pair $(\hat x,\hat X)$ to \eqref{maxWNP-relax}. It is defined as follows:
\[
gap(\hat x,\hat X) = \frac{| \hbox{tr} (P_0 \hat X) - \hat x^T P_0 \hat x |}{|\hat x^T P_0 \hat x|}.
\]
With regard to Theorem~\ref{prop-ekviv} a value of $gap(\hat x,\hat X)$ below the given small tolerance level indicates that $\hat x$ is indeed the optimal solution to the original problem \eqref{maxWN} and so the constructed function $\sigma\in \K^N$ is an optimal anisotropy function minimizing the anisotropic energy and satisfying quadratic constraints \eqref{second}. For the number of Fourier modes $N\approx 100$ the value of $gap(\hat x,\hat X)$ was less than $10^{-4}$.

\begin{table}
\renewcommand{\arraystretch}{1.3}
\caption{Dependence of the time complexity of computation with respect to the number of Fourier modes $N$ and its experimental order of time complexity (eotc) }
\label{tab1}

\scriptsize
\begin{center}
\begin{tabular}{c|cc}
\hline
$N$& CPU(s)&  eotc \\
\hline
50  &   5   &  --    \\
100 &  39   &  2.87  \\
150 & 124   &  2.86  \\
200 & 407   &  4.13  \\
250 &1609   &  4.32  \\
300 &2267   &  4.12  \\
\hline
\end{tabular}

\end{center}

\end{table}

\section{Conclusions}
In this paper we analyzed a novel method of enhanced semidefinite relaxation for solving a class of non-convex quadratic optimization problems. We applied this methodology to the practical problem of construction of the optimal anisotropy function minimizing the anisotropic energy for a given Jordan curve in the plane.

\section*{Acknowledgments}
The authors thank professors M. Halick\'a and M. Hamala and referees for their constructive comments and suggestions.

\ifCLASSOPTIONcaptionsoff
  \newpage
\fi

\end{document}